%% file: ALE_sinum_revision1.tex
\documentclass[hidelinks,onefignum,onetabnum]{siamart220329}
\usepackage{booktabs}

\renewcommand{\textcolor}[2]{#2}

\input{ex_shared}

\ifpdf
\hypersetup{
  pdftitle={ALE_sub1},
  pdfauthor={T.Lundquist, A. Malan, and J.Nordstr\"{o}m}
}
\fi

\begin{document}

\maketitle

% REQUIRED
\begin{abstract}
We present a novel framework based on semi-bounded spatial operators for analyzing and discretizing initial boundary value problems on moving and deforming domains.  This development extends an existing framework for well-posed problems and energy stable discretizations from stationary domains to the general case including arbitrary mesh motion. In particular, we show that an energy estimate derived in the physical coordinate system \textcolor{green}{is equivalent to} a semi-bounded property \textcolor{green}{with respect to} a stationary reference domain. The continuous analysis leading up to this result is based on a skew-symmetric splitting of the material time derivative, and thus relies on the property of integration-by-parts. Following this, a mimetic energy stable arbitrary Lagrangian-Eulerian framework for semi-discretization is formulated, based on approximating the material time derivative in a way consistent with discrete summation-by-parts.  Thanks to the semi-bounded property, a method-of-lines approach using standard explicit or implicit time integration schemes can be applied to march the system forward in time.  The same type of stability arguments applies as for the corresponding stationary domain problem, without regards to additional properties such as discrete geometric conservation.
As an additional bonus we demonstrate that \textcolor{blue}{discrete} geometric conservation, in the sense of exact free-stream preservation, can still be achieved \textcolor{blue}{in an automatic way} with the new framework. However, we stress that \textcolor{blue}{this is not necessary for stability.}
\end{abstract}

% REQUIRED

% REQUIRED

\section{Introduction}

Semi-discrete approximations, i.e.  where the time variable is left continuous,  are widely used in the design and stability analysis of numerical methods for partial differential equations (PDEs). This approach, which separates the approximation in space and time, is commonly known as the method-of-lines.  
For general linear PDE problems posed on stationary domains,  $L_2$ energy estimates of the solution can be related to well-posedness through the concept of semi-bounded spatial operators \cite{Gustafsson08, Kreiss_Lorenz04}.  For semi-discrete approximations, \textcolor{green}{the} analogous property of the system matrix allows for standard time integration methods (e.g. $A-$stable implicit or conditionally stable explicit Runge-Kutta methods) to be applied in a stable way using the method-of-lines.

In the case of explicit time stepping,  the choice of a stable time step size is traditionally based on either trial-and-error or a simplified spectral analysis,  even if the semi-discrete matrix is not diagonalizable.  
A more limited set of schemes can be shown to formally preserve the energy stability associated with arbitrary semi-bounded operators in space \cite{Levy_Tadmor1998, Ranocha18,Sun_Shu2019}.  However,  such stability results tend to be associated with serious time step limitations.  So-called strong stability preserving schemes \cite{Shu88,Gottlieb05} may also satisfy strict energy bounds,  however this relies on a sufficient amount of dissipation such that the Euler forward scheme is contractive.  Regardless of which particular method is applied, by using the method-of-lines approach one typically seeks in some way, either directly or indirectly, to take advantage of an underlying semi-bounded property of the continuous problem.
Despite the limitations discussed above,  the convenience and efficacy of the semi-discrete framework means that the method-of-lines has remained the most popular approach for solving PDE problems numerically over many decades. 

The situation is markedly different when it comes to PDE problems posed on moving and deforming spatial domains. The time dependent nature of the $L_2$ norm on such domains implies that energy estimates derived in the Eulerian (physical) coordinate system no longer correspond in an obvious way with a semi-bounded property of the spatial operator.  Consequently,  a naive application of the method-of-lines may not lead to stability even if it does so for the corresponding problem posed on a stationary domain.  Various different approaches to solve PDE problems numerically on moving domains have been developed, and are generally categorized under the rubric of arbitrary Lagrangian-Eulerian (ALE) methods,  signifying that the computational frame of reference is allowed to move in an arbitrary way following the motion of the mesh nodes.  \textcolor{blue}{The stability theory for ALE methods is in general less developed than for stationary domain methods. Moreover, even in those cases where formal stability proofs have been presented, only carefully selected implicit time marching methods such as low order one-step or multi-step methods in  \cite{Formaggia99,Farhat01,Formaggia04,Boffi04,Badia06} or high order space-time discretizations in \cite{Bonito13a,Bonito13b,nikkar15,Kopriva2016,Zhou2019,nikkar15, Kopriva2016} apply.}

%In the vast available literature on ALE methods however, with a few exceptions (see e.g. \cite{Formaggia99,Farhat01,Formaggia04,Boffi04,Badia06,Bonito13a,Bonito13b,Zhou2019,nikkar15, Kopriva2016}) there is a general lack of rigorous stability theory.  \textcolor{red}{Moreover,  even in those successful cases only carefully selected (implicit) time marching methods such as implicit Euler in \cite{Formaggia99} or space-time high order discretizations in \cite{nikkar15, Kopriva2016} can be used.}

%Perhaps most notable among these exceptions is that of space-time schemes using summation-by-parts (SBP) operators in both space and time \cite{nikkar15, Kopriva2016}.
%However, even in this case the analysis is fully discrete in nature, and only applicable to certain types of implicit time integration methods based on the application of SBP in time \cite{Nordstrom13,Lundquist14,boom15, Linders20_RK}.

\textcolor{red}{
In this paper we develop a new method-of-lines framework for analyzing and solving PDE problems numerically on moving and deforming domains.  Based on a skew-symmetric splitting of the material derivative, we show that an energy estimate derived in the physical coordinate system is equivalent to a semi-bounded property of the spatial operator when posed on a stationary reference domain. This introduction of a stationary reference domain is however only needed} for the purpose of analysis. In practice, no explicit knowledge of the associated mapping functions or metric terms is needed. \textcolor{red}{The semi-discrete framework for numerical solution similarly guarantees that an energy estimate derived on the semi-discrete level is (for linear problems) equivalent to a bounded system matrix spectrum. This allows for a method-of-lines approach using standard time integration methods}.

%In this paper we develop a new framework for analyzing and solving PDE problems numerically on moving and deforming domains.  It follows a method-of-lines approach by  \textcolor{green}{first} focusing on semi-discrete approximations of the material time derivative. We rely on the classical concept of semi-boundedness for linear well-posedness, and semi-discrete stability follows from this in a similar way as for stationary domains.  
%The analysis is fundamentally based on making a variable substitution such that the material time derivative assumes a skew-symmetric split form.
%Thanks to this, we can prove that energy estimates derived with respect to the physical coordinate system can be directly associated with semi-boundedness. The formulation of this result relies on the introduction of a stationary reference domain. However this is done strictly for the purpose of analysis. In practice,  no explicit knowledge of the associated mapping functions or metric terms is needed. 

Our approach \textcolor{green}{thus} supersedes the requirement of discrete geometric conservation,  i.e.  that constant free-stream solutions are preserved exactly over time \textcolor{blue}{by the fully discrete scheme.} This property has often been considered as a general prerequisite for stable and/or accurate ALE methods,  see e.g. the literature study in \cite{Etienne09}. \textcolor{blue}{Instead, we satisfy a semi-discrete version of the Reynolds transport theorem by discretizing the split form of the material time derivative using nodal summation-by-parts (SBP) operators. Although not necessary for stability in our framework,  we finally demonstrate that discrete geometric conservation (in above sense of free-stream preservation) can still be achieved in an automatic way following the application of general time marching schemes.}

%As part of the continuous analysis we demonstrate that the same material derivative splitting leads to a simple proof of the Reynolds transport theorem.  All the essential parts of the continuous analysis thus follow from applying integration-by-parts (IBP) in space. As a consequence of this, a mimetic ALE framework for semi-discrete stability, valid for many different types of numerical methods, can subsequently be formulated based on the corresponding discrete concept of summation-by-parts (SBP).
%The semi-discrete framework guarantees that a strictly non-growing solution norm is associated with a matrix spectrum strictly confined to the left half plane (following from the semi-bounded property), thus allowing for a method-of-lines approach using standard time integration methods. In particular, our approach supersedes the requirement of discrete geometric conservation,  i.e.  that constant free-stream solutions are preserved exactly over time. This property has often been considered as a general prerequisite for stable and/or accurate ALE methods,  see e.g. the literature study in \cite{Etienne09}.
%Although not necessary for stability in our framework,  we finally demonstrate that geometric conservation (free-stream preservation) can still be achieved in a constructive way using standard time integration methods.

The rest of this paper is organized as follows.  Section \ref{sec:cont} is dedicated to  demonstrating the equivalence between energy estimates derived in physical space on the one hand,  and semi-boundedness with respect to a stationary reference domain on the other.  In section \ref{sec:semidisc} we formulate a mimetic semi-discrete framework for stability which follows step by step the continuous analysis of section \ref{sec:cont}.
Next, we demonstrate the additional sometimes desirable (but not essential for stability) property of free-stream preservation in section \ref{sec:freestream}. \textcolor{red}{In section \ref{sec:num} we validate the stability and accuracy of the proposed method-of-lines approach through a numerical experiment,} and finally in section \ref{sec:concl} we draw conclusions.

\section{The continuous framework}\label{sec:cont}
For ease of notation we consider a scalar initial boundary value problem in $p$ space dimensions on the general form
\begin{equation}\label{eq:cont1} 
\begin{aligned}
u_t = \ & D(u,x,t)u+F(x,t), && 
 x  \in \Omega, && \textcolor{red}{t_0 < t < T}\\
B(u,x,t)u = \ & g(x,t), && x \in \partial\Omega , && \textcolor{red}{t_0 < t < T}\\
u(x,t_0) = \ & u_0(x),
\end{aligned}
\end{equation}
where $x = (x_1,\ x_2,\ \ldots, \ x_p)^T$ is the vector of physical coordinates,  $u_t$ denotes the partial time derivative of the dependent variable $u=u(x,t)$, and $D$ and $B$ denotes a differential operator (in general non-linear) and a boundary operator, respectively.  Moreover, $F$ and $g$ denote the solution independent forcing and boundary data, respectively. 
\textcolor{red}{Throughout the rest of the paper we will assume that the problem (\ref{eq:cont1}) is well-posed, with the solution $u$ belonging to the space $\mathcal{U}$ of sufficiently regular functions as required by the energy method.  In particular, for this purpose we define $\mathcal{U}$ as the set of all functions $\phi$ such that $\phi_t$, $\nabla \phi$ and $D\phi$ exist and are continuous on the closure of the computational domain (i.e. on $[t_0, T] \times\Omega\cup\partial\Omega$),  where $\nabla$ denotes the gradient operator in space.  }

Throughout this paper, we will make frequent use of inner product notation for both volume and surface integrals over $\Omega$. Given two $L_2$ integrable functions $\phi$ and $\psi$ over both $\Omega$ and $\partial \Omega$, we write
\begin{equation}\label{eq:L2_cont}
\big( \phi, \psi \big)_\Omega =\int_{\Omega} \phi\psi \ dV , \quad \big( \phi, \psi \big)_{\partial\Omega}  =\oint_{\partial \Omega}  \phi\psi \ dS.
\end{equation}
An indispensable tool for analyzing PDE problems with the energy method is the concept of integration-by-parts (IBP).  Let $n = (n_1, \ n_2,\ \ldots, \ n_p)^T$ denote the outward pointing unit normal to $\Omega$. The IBP formula for a single coordinate direction is then given by
\begin{equation}\label{eq:IBP_partial}
 \big(\phi, \psi_{x_i}\big)_{\Omega} = -  \big(\phi_{x_i}, \psi \big)_\Omega+ \big(\phi, n_i  \psi\big)_{\partial\Omega},
\end{equation}
where $\phi$ and $\psi$ are \textcolor{red}{continuously differentiable functions}, and $\phi_{x_i}$ and $\psi_{x_i}$ denote partial derivatives in the $x_i$-direction.

\subsection{Semi-boundedness: the stationary domain case}
Since, to the best of our knowledge,  semi-boundedness has not previously been considered in the context of moving meshes, we start by reviewing this classical PDE concept in the way that it is usually applied, i.e. for a stationary domain $\Omega$ in (\ref{eq:cont1}). In this case, a standard application of the energy method starts out from the central observation that
\begin{equation}\label{eq:en_cont0}
\frac{1}{2}\frac{d}{dt}\|u\|^2_\Omega = \big(u, u_t\big)_\Omega,
\end{equation}
simply due to the product rule.  For the purpose of linear well-posedness (and later in the  discrete setting, stability), we need only to consider the case with zero data, $F=0$,  $g=0$.  This naturally leads to the following definition.

\begin{definition}\label{def:semi-b_cont}
For the case of a stationary domain $\Omega$ in (\ref{eq:cont1}), we say that the spatial operator (i.e.  the combination of $D$ and $B$) is (maximally) semi-bounded if for all $\phi\in \mathcal{U}$ satisfying $B\phi = 0$,  the estimate
\begin{equation}\label{eq:semi-b_cont1}
\big(\phi, D\phi\big)_\Omega \leq \alpha \|\phi\|_{\Omega}^2
\end{equation}
holds for some constant $\alpha$ independent of $\phi$.  In addition, if $\alpha\leq 0$ we then say that the operator is strictly semi-bounded. 
\end{definition}
\begin{remark}\label{remark:maxi_semi-b}
The term maximally semi-bounded refers to the case where
a minimal number of boundary conditions are used to obtain the bound (\ref{eq:semi-b_cont1}). The notion of maximal semi-boundedness can be applied to show that the linearized problem is well-posed, see e.g. \cite{Kreiss_Lorenz04, Nordstr20_nbc} for more details.  In the rest of the paper we will always assume that no overspecification is done, and for ease of notation we simply refer to such problems as being semi-bounded.
\end{remark}

Note that (\ref{eq:semi-b_cont1}) inserted into (\ref{eq:en_cont0}) yields, for $F=g=0$,
\begin{equation}\label{eq:en_cont_stat}
\frac{1}{2}\frac{d}{dt}\|u\|^2_\Omega \leq \alpha \|u\|_{\Omega}^2.
\end{equation}
Thus, for a stationary domain there is a direct
correspondence between energy estimates of the solution and semi-boundedness of the operator according to Definition \ref{def:semi-b_cont}.

For semi-discrete approximations of the type $\boldsymbol{U}_t = \mathcal{M}\boldsymbol{U}$, note that an analogous property to (\ref{eq:semi-b_cont1}) is given by the matrix estimate
\begin{equation}\label{eq:semi-b_disc1}
\mathcal{H}\mathcal{M} + \mathcal{M}^T\mathcal{H} - \textcolor{red}{2\alpha \mathcal{H} \leq 0,} 
\end{equation}
leading to 
\begin{equation}\label{eq:en_disc_stat}
\frac{1}{2}\frac{d}{dt}\|\boldsymbol{U}\|^2_{\mathcal{H}} \leq \alpha \|\boldsymbol{U}\|_{\mathcal{H}}^2,
\end{equation}
where \textcolor{red}{$\mathcal{H}=\mathcal{H}^T$ is a positive definite matrix defining a \textcolor{red}{discrete equivalent to the} $L_2$ inner product,} $( \boldsymbol{\Phi}, \boldsymbol{\Psi})_\mathcal{H} =\boldsymbol{\Phi}^T\mathcal{H} \boldsymbol{\Psi}$. In particular, if $\alpha \leq 0$ (i.e. strict semi-boundedness) then condition (\ref{eq:semi-b_disc1}) implies that the spectrum of $\mathcal{M}$ is confined to the left half plane, allowing for the stable application of standard time integration methods. 

If the domain is not stationary, i.e. $\Omega = \Omega(t)$, then (\ref{eq:en_cont0}) no longer holds,  and thus the direct correspondence between semi-boundedness (\ref{eq:semi-b_cont1}) and energy estimates of the form (\ref{eq:en_cont_stat}) is lost.  Similarly for the semi-discrete case, $\mathcal{H}=\mathcal{H}(t)$ invalidates the direct correspondence between (\ref{eq:semi-b_disc1}) and (\ref{eq:en_disc_stat}). In other words, there is no longer a clear connection between semi-discrete and fully discrete stability using the method-of-lines.
Our main goal will therefore be to reestablish the direct link between energy estimates and the concept of semi-boundedness for the general case with a moving domain, $\Omega = \Omega(t)$. 
After semi-discretization,  this approach will lead to schemes satisfying (\ref{eq:semi-b_disc1}) by construction, thus allowing for a stable application of the method-of-lines (again assuming $\alpha \leq 0$).

\subsection{Domain motion}\label{sec:motion}
Next, consider the general case $\Omega = \Omega(t)$ in (\ref{eq:cont1}), i.e. the domain boundary $\partial \Omega$ is subject to an arbitrary motion over time.  \textcolor{green}{For simplicity} we \textcolor{green}{will} assume that this motion is given as data to the problem, and thus independent of the solution $u$. 
When solving (\ref{eq:cont1}) numerically on a mesh that moves together with $\Omega$, it is natural to study solutions in the reference frame of the moving mesh nodes rather than in the original, physical coordinates $x$.  For this purpose we formally introduce a reference coordinate system wherein the mesh nodes are stationary. In terms of these new coordinates we denote the spatial domain with $\hat{x}\in\hat{\Omega}$, where the reference domain $\hat{\Omega}$ is stationary in time.

We impose no restrictions on the mesh motion apart from non-degeneracy in the sense that
%, for each time $t$, 
the existence of a bijective  
\textcolor{red}{and sufficiently smooth} 
%\textcolor{red}{and piecewise twice continuously differentiable} 
coordinate mapping $x=x(\hat{x},t)$, $\hat{x}=\hat{x}(x,t)$ is assumed.  
\textcolor{red}{Assuming further that the computational domain is subdivided into a finite number of  non-overlapping blocks or elements,  we require that this mapping is globally continuous as well as twice continuously differentiable on the closure of each such subdomain.}
The evolution of the physical domain $\Omega$ over time can \textcolor{red}{then} be described  by the \textcolor{red}{piecewise smooth} coordinate velocity $\dot{ x}=(\dot{x}_1,\ \dot{x}_2, \ \ldots, \ \dot{x}_p)^T$, defined by 
\begin{equation}\label{eq:coor_velo}
\dot{x}_i=\dot{x}_i(x, t) =: \frac{\partial }{\partial t}x_i(\textcolor{red}{\hat{x}(x,t)}, t).
\end{equation}
The functions $\dot{ x}_i(x, t)$ in (\ref{eq:coor_velo}) \textcolor{green}{will henceforth be} assumed to represent known data describing the mesh motion, while no explicit knowledge of the \textcolor{green}{underlying} mapping functions $x=x(\hat{x},t)$, $\hat{x}=\hat{x}(x,t)$ or any spatial derivatives thereof (metric terms) will be required.
\begin{remark}\label{remark:curvi}
In ALE formulations it is common to employ the initial domain as reference, i.e.  $\hat{\Omega}  = \Omega(t_0)$ and $x(\hat{x},t_0)=\hat{x}$ \cite{Belytschko2000}.  
While this choice is always possible to make, it is not necessary. For example, if curvilinear coordinates are employed to describe $\Omega$, it may be more convenient to make use of this known transformation instead, and thus define $\hat{\Omega}$ \textcolor{red}{piecewise on each computational subdomain in terms of the unit hypercube}.
\end{remark}

\subsection{The material time derivative}\label{sec:material_cont}
Note that the (arbitrary) reference domain $\hat{\Omega}$ is the only stationary domain we explicitly consider in this work, and thus there is no need to make a distinction between $\hat{\Omega}$ and the alternative Lagrangian (material) view of \textcolor{green}{particle motion} in continuum mechanics.  Without risk of ambiguity, we will therefore refer to \textcolor{green}{all} time derivatives evaluated in the reference coordinate system (i.e. following the moving mesh nodes) as material time derivatives.

In terms of the reference coordinates $\hat{x}$, the so-called material (or total) time derivative $d\phi/dt$ of a general function $\phi$ can be obtained by application of the chain rule, i.e. we have
\begin{equation}\label{eq:material_cont1}
\frac{d\phi}{dt} =\frac{d }{d t} \phi (x(\hat{x},t), t)=  \phi_t +\dot{ x}^T
\big( \nabla \phi \big)  ,
\end{equation}
where $\dot{x}$ is defined in (\ref{eq:coor_velo})\textcolor{green}{, and} $\phi_t = \partial \phi (x, t)/ \partial t $ 
as before denotes a partial time derivative with respect to the physical coordinates $x$. 
By substituting $u_t$ for $du/dt$  in (\ref{eq:cont1}), the mesh motion can thus be seen to induce an additional spatial term into the PDE.

We start by considering a splitting \cite{ref:NORD062} of the spatial term in (\ref{eq:material_cont1}) into symmetric and skew-symmetric parts.  Applying the product rule, we get
\begin{equation}\label{eq:material_cont_split}
\frac{d\phi}{dt} =  \phi_t + 
D_m\phi   - \frac{1}{2} \big(\nabla\cdot \dot{ x} \big) \phi,
\end{equation}
where the operator $D_m$ is defined by
\begin{equation}\label{eq:Q_cont}
D_m\phi = \frac{1}{2}\Big[\dot{ x}^T \big( \nabla \phi \big)   +\nabla\cdot \big(\dot{ x}  \phi \big)  \Big] .
\end{equation}
From the partial derivative formulae (\ref{eq:IBP_partial}), the combined IBP property
\begin{equation}\label{eq:IBP_material}
\big( \phi,D_m\psi\big)_\Omega = -\big(  D_m\phi,\psi  \big)_\Omega +
 \big(\phi, \big(n^T\dot{ x}\big) \psi \big)_{\partial \Omega}
\end{equation}
now follows, showing that $D_m$ is a skew-symmetric operator in the interior of $\Omega$.  

The divergence $\nabla\cdot \dot{ x}$ can be both positive or negative depending on whether local volumes are increasing or decreasing, and the last term in (\ref{eq:material_cont_split}) thus constitutes a source term with an indefinite sign. A naive semi-discretization of (\ref{eq:cont1}) will include the same indefinite contribution from the mesh motion, causing even $A-$stable time integration methods to fail.  
In order to resolve this problem, our goal will be to establish a general link between the \textcolor{green}{two} concepts of energy stability and semi-boundedness,  \textcolor{green}{which to our knowledge has} previously only been considered for stationary domains.

\begin{remark}
For the special case of volume preserving (isochoric) motion, the source term involving $\nabla\cdot \dot{ x}$ is always zero, making this case significantly easier to analyze and discretize, see \cite{Lundquist20_JCP}.
\end{remark}

\subsection{Variable substitution}

By the inverse function theorem of multivariate calculus, \textcolor{green}{ for the mesh motion to be non-degenerate (i.e. bijective)} it is sufficient that the forward mapping function $x=x(\hat{x}, t)$ has a strictly positive Jacobian determinant,
\begin{equation}
\label{eq:inv_fun_thm_cont}
J(x,t) =\mathrm{det}\Bigg(\frac{\partial x}{\partial \hat{x}}\Bigg) >  0.
\end{equation}
\textcolor{red}{In fact, the smoothness assumption on $x(\hat{x}, t)$ further implies that $J$ must be uniformly bounded from both below and above}.
This function $J$ will play a central role in the upcoming analysis of (\ref{eq:cont1}) for the general case $\Omega  = \Omega (t)$.
Recall that we consider the coordinate velocity vector $\dot{x}(x,t)$ as data to the problem,  while no explicit knowledge of the Jacobian matrix $\partial x /\partial \hat{x}$ itself is assumed.  
It is therefore important to note that the value of $J$ in (\ref{eq:inv_fun_thm_cont}) can be directly related to $\dot{x}$ through the initial value problem
\begin{equation}\label{eq:GCL_cont} 
\begin{aligned}
\frac{dJ}{dt} = \ &\big(\nabla\cdot\dot{ x}\big)J, && 
 x  \in \Omega (t), && t\geq t_0 \\
J(x,t_0) = \ & J_0(x) >0, && 
 x  \in \Omega (t).
\end{aligned}
\end{equation}
In other words, local volumes change at a rate which is proportional to the divergence of the coordinate velocity vector.  
%\textcolor{red}{
%Note that (\ref{eq:GCL_cont}) holds as long as the mapping is sufficiently regular in order for the Jacobian matrix $\partial x /\partial \hat{x}$ to be pieciewise differentiable in time.  
A proof of this result can be found e.g. in \cite{gurtin81}, page 77.
%}

Next, we demonstrate that the function $J$ in (\ref{eq:GCL_cont}) can be employed in order to eliminate the indefinite source term previously identified in the material derivative splitting (\ref{eq:material_cont_split}).
For a general function $\phi$ defined on $\Omega (t)$, consider the corresponding function $\hat{\phi}$ \textcolor{red}{on the reference domain $\hat{\Omega}$ given by the variable substitution} 
\begin{equation}\label{eq:varsub}
\textcolor{red}{\hat{\phi}(\hat{x},t) = \sqrt{J(x(\hat{x},t),t)}\phi(x(\hat{x}, t), t).}
\end{equation}
This leads to
\begin{equation*}
\frac{d\hat{\phi} }{dt} =  \sqrt{J}\frac{d\phi }{dt} +\frac{d\sqrt{J}}{dt} \phi .
\end{equation*}
Moreover, since $J>0$ we can rewrite (\ref{eq:GCL_cont}) into the equivalent form
\begin{equation*} 
\frac{d \sqrt{J}}{dt} =  \frac{1}{2} \big(\nabla\cdot\dot{ x}\big)
\sqrt{J}.
\end{equation*}
After also inserting (\ref{eq:material_cont_split}), we thus get
\begin{equation*}
\frac{d\hat{\phi}}{dt} =
\sqrt{J}\Big[\phi_t+D_m\phi   - \frac{1}{2} \big(\nabla\cdot \dot{ x} \big) \phi\Big]+\frac{1}{2} \big(\nabla\cdot\dot{ x}\big)
\sqrt{J} \phi.
\end{equation*} 
Note that the same indefinite source term $\frac{1}{2} \big(\nabla\cdot\dot{ x}\big)
\sqrt{J} \phi$ now appears in two places with opposite signs on the right hand side above, and thus cancel. We have proven
\begin{lemma}\label{lemma:materialderiv}
Let the domain $\Omega(t)$ be subject to a non-degenerate (\ref{eq:inv_fun_thm_cont}) motion described by the \textcolor{red}{piecewise smooth coordinate velocity function} $\dot{x}$ in (\ref{eq:coor_velo}).
For any differentiable function $\phi$ on $\Omega$, the corresponding function $\hat{\phi}= \sqrt{J}\phi$ \textcolor{red}{on $\hat{\Omega}$ (see (\ref{eq:varsub}))} then satisfies
\begin{equation}\label{eq:material_cont_modif}
\frac{d\hat{\phi}}{dt} =
\sqrt{J}\big(\phi_t+D_m\phi
\big),
\end{equation}
where $D_m$ is the skew-symmetric operator defined in (\ref{eq:Q_cont}), and $J>0$ satisfies (\ref{eq:GCL_cont}).
\end{lemma}

Going back to the original model problem (\ref{eq:cont1}),  the variable substitution \textcolor{red}{(\ref{eq:varsub}) with $\phi=u$} now yields, using Lemma \ref{lemma:materialderiv},
\begin{equation*}
\frac{d\hat{u}}{dt} =
\sqrt{J}\big(Du + F+D_m u
\big).
\end{equation*}
Stated in terms of the stationary reference domain $\hat{\Omega}$, the problem (\ref{eq:cont1}) can thus be equivalently written as
\begin{equation}\label{eq:cont_MOL} 
\begin{aligned}
\frac{d\hat{u}}{dt} = \ &  \hat{D}(\hat{u}, \hat{x}, t)\hat{u} + \sqrt{J}F(x(\hat{x},t), t),  && 
 \hat{x}  \in \hat{\Omega}, && t_0<t<T \\
\hat{B}(\hat{u},  \hat{x}, t) \hat{u} = \ & g(x(\hat{x},t), t)  && \hat{x} \in \partial\hat{\Omega} , && t_0<t<T\\
\hat{u}(\hat{x},t_0) = \ & \sqrt{J}u_0(x(\hat{x},t), t),&& \hat{x}  \in \hat{\Omega}, 
\end{aligned}
\end{equation}
where
\begin{equation*}
\hat{D}=   \sqrt{J} \big(D +D_m\big)\sqrt{J} ^{-1}, \quad \hat{B}=  B\sqrt{J}^{-1}.
\end{equation*}
%\textcolor{red}{Note that since $\dot{x}$ is only assumed to be piecewise smooth, even if $u\in \mathcal{U}$, the solution $\hat{u}=\sqrt{J}u$ to (\ref{eq:cont_MOL}) may still only exist in a weak sense.}
\textcolor{red}{Moreover, if $u\in\mathcal{U}$ holds then the solution $\hat{u}$ in (\ref{eq:cont_MOL}) belongs to the space $\hat{\mathcal{U}}$ given by all $\hat{\phi}$ in (\ref{eq:varsub}) such that $\phi\in\mathcal{U}$. 
Note that since the smoothness properties of the coordinate mapping is only assumed to hold piecewise, the functions contained in $\hat{\mathcal{U}}$ can be locally discontinuous at block or element boundaries.  However, the operator $\hat{D}$ is still well-defined for all $\hat{u}\in\hat{\mathcal{U}}$.
Next,} we will proceed to demonstrate that an energy estimate of the type (\ref{eq:en_cont_stat}) with respect to the moving domain $\Omega(t)$ \textcolor{green}{is equivalent to the} semi-boundedness (\ref{eq:semi-b_cont1}) \textcolor{green}{of $\hat{D}$} on $\hat{\Omega}$. To do this, we will also need the Reynolds transport theorem.

\subsection{The Reynolds transport theorem}\label{sec:Reynolds_cont}
For the purpose of applying the energy method to (\ref{eq:cont1}),  an analogous property to (\ref{eq:en_cont0}) for the case $\Omega = \Omega (t)$ can be obtained from the well-known Reynolds transport theorem (see e.g.  \cite{Leal2007}).
In the interest of later formulating a semi-discrete analogue to this result, we give a simple (and to our knowledge, novel) proof of the continuous Reynolds transport theorem below based on the application of Lemma \ref{lemma:materialderiv}.  \textcolor{blue}{This derivation will later serve as a guide for the discrete analysis in section \ref{sec:Reynolds_disc}.} For two functions $\phi$ and $\psi$, we first note that by defining $\hat{\phi}= \sqrt{J}\phi$ and $\hat{\psi}= \sqrt{J}\psi$ as in (\ref{eq:varsub}), the identity
 \begin{equation}\label{eq:norm_equi_cont}
 \big(\phi,   \psi\big)_\Omega =   \big(\phi,  J \psi\big)_{\hat{\Omega}}  = \big(\hat{\phi}, \hat{\psi}\big)_{\hat{\Omega}}
 \end{equation}
holds.
Since $\hat{\Omega}$ is stationary in time, we can \textcolor{green}{thus} write the time derivative of a general inner product on $\Omega$ as
 \begin{equation*}
 \frac{d}{d t}\big(\phi, \psi\big)_{\Omega} = \frac{d}{d t}\big(\hat{\phi}, \hat{\psi}\big)_{\hat{\Omega}}  = \Big( \frac{d\hat{\phi}}{d t}, \hat{\psi}\Big)_{\hat{\Omega}} + \Big( \hat{\phi}, \frac{d\hat{\psi}}{d t}\Big)_{\hat{\Omega}} .
 \end{equation*}
Now inserting (\ref{eq:material_cont_modif}) from Lemma \ref{lemma:materialderiv}, the first term on the right hand side above becomes
 \begin{equation*}
 \Big( \frac{d\hat{\phi}}{d t}, \hat{\psi}\Big)_{\hat{\Omega}} =  \Big( \sqrt{J}\phi_t, \hat{\psi}\Big)_{\hat{\Omega}} +  \Big( \sqrt{J}D_m\phi, \hat{\psi}\Big)_{\hat{\Omega}} =  \Big( \phi_t, \psi\Big)_\Omega+  \Big( D_m\phi, \psi\Big)_\Omega,
 \end{equation*}
and the second term can of course be written in a similar way.  From this and the IBP property (\ref{eq:IBP_material}), the Reynolds transport theorem
 \begin{equation}\label{eq:Reynolds2_cont}
\frac{d}{dt} \big( \phi, \psi \big)_\Omega = 
\big(\phi_t, \psi \big)_\Omega +
\big( \phi, \psi_t\big)_\Omega +
\big( \phi,\big(n^T\dot{ x}\big) \psi \big)_{\partial \Omega }
\end{equation}
follows, which we will use in order to apply the energy method to the original formulation of the PDE problem in (\ref{eq:cont1}). 

\subsection{The energy method}\label{sec:analysis}

By inserting $\phi = \psi = u$ into the Reynolds transport theorem (\ref{eq:Reynolds2_cont}), we get
\begin{equation}\label{eq:en_cont}
\frac{d}{dt} \|u\|_\Omega^2
=2 \big( u, u_t \big)_\Omega +
\big( u,\big(n^T\dot{ x}\big)   u \big)_{\partial\Omega},
\end{equation}
corresponding to (\ref{eq:en_cont0}) in the stationary domain case.
Inserting the PDE in (\ref{eq:cont1}) into (\ref{eq:en_cont}) now yields
\begin{equation*}
\frac{d}{dt}\|u\|^2_\Omega =2 \big(u, Du\big)_\Omega+2\big(u, F)_\Omega+\big( u, \big(n^T\dot{ x}\big)   u \big)_{\partial\Omega}.
\end{equation*}
In order to arrive at the same type of energy estimate as before in (\ref{eq:en_cont_stat}) for the zero data case $F = g = 0$, the combination of $D$ and $B$ must clearly satisfy the following property:
For all \textcolor{red}{functions $\phi\in\mathcal{U}$} satisfying $B\phi = 0$,  the estimate
\begin{equation}\label{eq:enboun_cont1}
\big(\phi, D\phi\big)_\Omega+\frac{1}{2}\big( \phi, \big(n^T\dot{ x}\big)   \phi \big)_{\partial\Omega} \leq \alpha\|\phi\|^2_\Omega
\end{equation}
must hold. Note that this immediately implies the energy estimate (\ref{eq:en_cont_stat}) for the case $F = g = 0$. We will therefore in general assume that the PDE model including boundary condition can be designed in such a way that (\ref{eq:enboun_cont1}) holds.
 However, since $n^T\dot{ x}$ can have both positive or negative sign,  (\ref{eq:enboun_cont1}) does \textit{not} characterize the spatial operator as semi-bounded with respect to $\alpha$ according to Definition \ref{def:semi-b_cont}.
As already discussed, semi-boundedness (in particular with $\alpha\leq 0$) is a necessary property for stable time integration with the method-of-lines.  In its current form, the estimate (\ref{eq:enboun_cont1}) is therefore insufficient for this purpose.
We are now ready to state the first main result of this paper.
\begin{proposition}\label{thm:semi-b}
For a moving domain $\Omega(t)$ the following two statements are equivalent:
\begin{enumerate}
 \item The \textcolor{green}{energy} estimate (\ref{eq:enboun_cont1}) holds for all \textcolor{red}{$\phi\in \mathcal{U}$} satisfying $B\phi = 0$.
\item The spatial operator in (\ref{eq:cont_MOL}) \textcolor{red}{satisfies the semi-bounded property}
\begin{equation}\label{eq:semi-b_cont2}
\big(\hat{\phi}, \hat{D}\hat{\phi}\big)_{\hat{\Omega}} \leq \alpha \|\hat{\phi}\|^2_{\hat{\Omega}}
\end{equation}
%\textcolor{red}{for all $\phi\in \mathcal{U}$ satisfying $B\phi = 0$, where $\hat{\phi}=\sqrt{J}\phi$}.
 \textcolor{red}{for all $\hat{\phi}\in\hat{\mathcal{U}}$ satisfying $\hat{B}\hat{\phi} = 0$. }
\end{enumerate}
\begin{proof}
First note that \textcolor{red}{for any $\hat{\phi}$ defined as in (\ref{eq:varsub}), the condition} $B\phi = 0$ is equivalent to $\hat{B}\hat{\phi} = 0$. 
Moreover, \textcolor{red}{by definition of $\hat{D}$ in (\ref{eq:cont_MOL}) we have}
\begin{equation*}
\big(\hat{\phi}, \hat{D}\hat{\phi}\big)_{\hat{\Omega}}  = \big(\sqrt{J}\phi, \sqrt{J}\big(D \phi+D_m\phi\big)\big)_{\hat{\Omega}} = \big(\phi, D \phi+D_m\phi\big)_\Omega.
\end{equation*}
Applying the IBP property (\ref{eq:IBP_material}), this further yields
\begin{equation*}
\big(\hat{\phi}, \hat{D}\hat{\phi}\big)_{\hat{\Omega}}  = \big(\phi, D\phi\big)_\Omega+\frac{1}{2}\big( \phi, \big(n^T\dot{ x}\big)   \phi \big)_{\partial\Omega}.
\end{equation*}
We have thus shown that the left hand sides of (\ref{eq:enboun_cont1}) and (\ref{eq:semi-b_cont2}) are identical. Since the right hand sides are also the same by (\ref{eq:norm_equi_cont}),  the two conditions are indeed equivalent.
\end{proof}
\end{proposition}

\textcolor{green}{In light of this result}, the energy analysis for any specific PDE problem can be carried out directly with respect to the original formulation (\ref{eq:cont1}) in physical space. Once an estimate of the form (\ref{eq:enboun_cont1}) has been derived,  semi-boundedness of the equivalent reference domain formulation (\ref{eq:cont_MOL}) then follows \textcolor{green}{automatically} by Proposition \ref{thm:semi-b}.
For semi-discrete approximations, recall that the analogous property to semi-boundedness given by (\ref{eq:semi-b_disc1}) implies a matrix spectrum confined to the left half plane if $\alpha \leq 0$.  The main goal of the semi-discrete ALE framework presented in section \ref{sec:semidisc} will therefore be to establish a similar equivalent relationship between a discrete energy estimate and the matrix condition (\ref{eq:semi-b_disc1}).

\subsection{A note on boundary condition design and lifting operators}
In order to simplify the design and analysis of boundary conditions leading to an estimate of the form (\ref{eq:enboun_cont1}), a so-called lifting operator technique \cite{Nordstrom17, Arnold02} may be employed. 
This approach also has the benefit of later providing a direct path to implementing stable boundary conditions in the semi-discrete scheme.

Thus,  consider including the boundary condition in (\ref{eq:cont1}) directly on the right hand side of the PDE using the penalty formulation
\begin{equation}\label{eq:cont1_lifting} 
\begin{aligned}
u_t = \ & Du+L_\delta(Bu-g)+F, && 
 x  \in \Omega (t), && \textcolor{red}{t_0<t< T} \\
u(x,t_0) = \ & u_0(x).
\end{aligned}
\end{equation}
The so-called lifting operator $L_\delta$ is here defined with respect to some operator $ \delta$ such that for any $\phi$, $\psi$, we have
\begin{equation}\label{eq:lifting_cont}
\big(\phi,   L_\delta \psi\big)_\Omega = \big(\delta \phi,  \psi\big)_{\partial{\Omega}}.
\end{equation}
Note that the term scaled by $L_\delta$ in (\ref{eq:cont1_lifting}) is here only added for the purpose of analysis; \textcolor{green}{as long as} the boundary condition \textcolor{green}{holds with equality}, then this added term has the value zero.
Applying the energy method to (\ref{eq:cont1_lifting}) directly (with zero data $g=0$, $F=0$) now yields
\begin{equation*}
\frac{d}{dt} \|u\|_\Omega^2 = \mathcal{E}(D,B,u),
\end{equation*}
where the function $\mathcal{E}$ is defined as
\begin{equation}\label{eq:eps_cont}
\mathcal{E}(D,B,\phi) =\big(\phi, D\phi\big)_\Omega+\frac{1}{2} \big(  \phi, \big(n^T\dot{ x}\big)   \phi \big)_{\partial\Omega}+\big(\delta\phi, B\phi\big)_{\partial\Omega} .
\end{equation}
Note in particular that if we insert $B\phi = 0$, then $\mathcal{E}(D,B,\phi)$ exactly reduces to the left hand side of (\ref{eq:enboun_cont1}).  Thus, \textcolor{red}{if there exists an operator $\delta$ such that} 
\begin{equation}\label{eq:enboun_cont2}
\mathcal{E}(D,B,\phi) \leq \alpha \|\phi\|_\Omega^2
\end{equation}
holds for all \textcolor{red}{$\phi\in\mathcal{U}$}, then this is clearly a sufficient condition for (\ref{eq:enboun_cont1}) to hold for all \textcolor{red}{$\phi\in\mathcal{U}$}, satisfying $B\phi=0$.  By combining the differential and boundary operators $D$ and $B$ into a single expression, the weak condition (\ref{eq:enboun_cont2}) can be \textcolor{green}{easier} to evaluate in practice than the original strong condition (\ref{eq:enboun_cont1}).

\subsubsection{Example}\label{sec:example}
To illustrate the constructive use of \textcolor{green}{the weak condition} (\ref{eq:enboun_cont2}) for boundary condition design with an example, we consider the constant coefficient advection-diffusion problem on a half-line, given by
 \begin{equation*}
 \begin{aligned}
 u_t  = \ & \epsilon u_{xx} - au_x, && x < x_e(t),  &&t>0 \\
Bu = \ & g(t), && x =x_e(t), && t> 0 \\
u(x, 0) = \ &u_0(x),  && x < x_e(t), 
 \end{aligned}
 \end{equation*}
where the type of boundary condition will be determined by the choice of operator $B$.
For this PDE problem,  $\mathcal{E}$ in (\ref{eq:eps_cont}) becomes
 \begin{equation*}
 \begin{aligned}
\mathcal{E} = \ & \int^{x_e} \phi(\epsilon \phi_{xx}-a\phi_x) dx +\Big[\big(\delta \phi\big)\big( B \phi\big)+ \frac{1}{2}\dot{x}_e \phi^2 \Big]_{x=x_e} \\
=  \ & -\epsilon \|\phi_x\|^2 + \Big[\big(\delta \phi\big)\big( B \phi\big)+\epsilon \phi\phi_x+\frac{1}{2}\big(\dot{x}_e-a\big)\phi^2 \Big]_{x=x_e},
\end{aligned}
\end{equation*}
and we recall that $\delta$ is some operator still to be specified. 

For example, by imposing the Dirichlet boundary condition $u=g(t)$ at $x=x_e$, we can specify
\begin{equation*}
B = 1 , \quad \delta = \frac{1}{2}(a-\dot{x}_e)-\epsilon\partial_x \quad \Rightarrow \quad \mathcal{E} = -\epsilon \|\phi_x\|^2,
\end{equation*}
thus satisfying (\ref{eq:enboun_cont2}) with $\alpha = 0$.  Recall that this is a sufficient condition for (\ref{eq:enboun_cont1}) to hold with $B\phi = 0$. Hence, Proposition \ref{thm:semi-b} yields the conclusion that also (\ref{eq:semi-b_cont2}) is satisfied for the same value $\alpha  = 0$, i.e.  we obtain  strict semi-boundedness.  Alternatively, by imposing the characteristic boundary condition $ \frac{1}{2}\big( a-\dot{x}_e -|a-\dot{x}_e|\big)u-\epsilon \phi_x = g(t)$, we can set
\begin{equation*}
B = \frac{1}{2}\big( a-\dot{x}_e -|a-\dot{x}_e|\big)-\epsilon  \partial_x, \quad \delta = 1 \quad \Rightarrow \quad \mathcal{E} =-\epsilon \|\phi_x\|^2 - \textcolor{green}{\frac{1}{2}}|a-\dot{x}_e|\phi(x_e,t)^2,
\end{equation*}
thus yielding the same conclusion as for the Dirichlet condition above. 

Furthermore, notice that in the case of pure advection,  i.e. $\epsilon = 0$, the characteristic boundary operator $B$ above is designed to only be non-zero at times when $\dot{x}_e - a$ yields a positive contribution to $\mathcal{E}$, i.e. when $x_e$ represents an inflow boundary. This shows that the corresponding operator in (\ref{eq:cont_MOL}) is maximally semi-bounded, and the problem is thus well-posed, see Remark \ref{remark:maxi_semi-b}.

\section{The semi-discrete framework}\label{sec:semidisc}
Recall that the result of Lemma \ref{lemma:materialderiv} in the continuous analysis followed from the velocity based definition (\ref{eq:GCL_cont}) of the Jacobian determinant $J$. Moreover, the proof of the Reynolds transport theorem given in section \ref{sec:Reynolds_cont} further relied on the application of \textcolor{green}{IBP (\ref{eq:IBP_partial})}. In order to mimick the analysis leading up to Proposition \ref{thm:semi-b} for semi-discrete approximations,  we need to focus our attention on these two properties. 
%\textcolor{green}{Hence, as} long as (\ref{eq:IBP_partial}) and (\ref{eq:GCL_cont}) hold discretely, then the semi-discrete analysis can be carried out in an identical way \textcolor{green}{for all such discretizations schemes.}
\textcolor{green}{Hence, the semi-discrete analysis will be carried out in an identical way for all types of numerical schemes that satisfy (\ref{eq:IBP_partial}) and (\ref{eq:GCL_cont}) discretely.}

For discretization, we consider a volume mesh with $n_V$ nodes representing $\Omega$, such that a subset of these nodes also defines a surface mesh of size $n_S$ representing $\partial \Omega$. We also define a restriction operator $E$ from the former mesh to the latter. That is, if $\Phi$ is the restriction of some vector quantity $\boldsymbol{\Phi}$ from all nodes to only the surface nodes, we then express this relationship as\begin{equation*}
\Phi = E \boldsymbol{\Phi},
\end{equation*}
where $E$ is a $n_S\times n_V$ matrix in which each row only contains a single non-zero element (with the value $1$). As a general convention,  we employ underscore notation to indicate the diagonal matrix obtained by inserting the underlined vector along the diagonal. That is, for vectors $\boldsymbol{\Phi}$ and $\Phi$ of dimension $n_V$ and $n_S$ respectively,  we frequently use the notation
\begin{equation*}
\nwu{\boldsymbol{\Phi}} =: \mathrm{Diag}(\boldsymbol{\Phi}), \quad \nwu{\Phi} =: \mathrm{Diag}(\Phi).
\end{equation*}

For an arbitrary set of vectors $\boldsymbol{\Phi}$, $\boldsymbol{\Psi}$, $\Phi$ and $\Psi$ with consistent dimensions, we consider the following discrete inner product definitions corresponding to (\ref{eq:L2_cont}) in the continuous analysis
\begin{equation}\label{eq:L2_disc_alt}
\big( \boldsymbol{\Phi}, \boldsymbol{\Psi} \big)_{{\mathcal{P}}} = \sum_{j=1}^{n_V} \boldsymbol{\Phi}_j\boldsymbol{\Psi}_j{\mathcal{P}}_j=\boldsymbol{\Phi}^T\nwu{\mathcal{P}}\boldsymbol{\Psi} , \quad 
\big( \Phi, \Psi \big)_P = \sum_{i=1}^{n_S} \Phi_j\Psi_j P_j,=\Phi^T\nwu{P}\Psi  ,
\end{equation}
where $\mathcal{P}=(\mathcal{P}_1, \mathcal{P}_2, \ldots, \mathcal{P}_{n_V})^T$ and $P=(P_1, P_2, \ldots, P_{n_S})^T$  denote two sets of quadrature weights for the volume and surface mesh, respectively.

\subsection{SBP partial derivative operators in space}
Since we aim to mimic the continuous analysis, we first need a discrete equivalent to IBP, referred to as summation-by-parts (SBP).  We thus consider a set of linear operators $\mathcal{D}_{x_i}$ approximating $\frac{\partial }{\partial x_i}$ on the mesh nodes,  satisfying the following SBP properties corresponding to (\ref{eq:IBP_partial})
\begin{equation}\label{eq:SBP_partial0}
  \big(\boldsymbol{\Phi}, \mathcal{D}_{x_i} \boldsymbol{\Psi} \big)_{\mathcal{P}} = -  \big(\mathcal{D}_{x_i}\boldsymbol{\Phi}, \boldsymbol{\Psi} \big)_{\mathcal{P}}+ \big(\Phi, \nwu{N}_{i}  \Psi \big)_P,
\end{equation}
\textcolor{green}{where $\Phi = E \boldsymbol{\Phi}$, $\Psi = E \boldsymbol{\Psi}$, and}
where $N_{i}$ contains the surface node values of the normal vector component $n_{i}$ (either exact or approximated). \textcolor{red}{Note that we do not make any assumptions regarding the function space of admissible numerical solutions. Rather, (\ref{eq:SBP_partial0}) should hold for arbitrary vectors $\boldsymbol{\Phi} $, $\boldsymbol{\Psi} $ of size $n_V$. }

Indeed,  \textcolor{green}{properties of the type (\ref{eq:SBP_partial0})} can be satisfied by design for a wide range of \textcolor{green}{nodal} numerical \textcolor{green}{schemes}, employing both unstructured as well as structured curvilinear meshes. The specific quadratures ${\mathcal{P}}$ and $P$ appearing in (\ref{eq:SBP_partial0}) are \textcolor{blue}{in general uniquely} determined by the particular scheme. Examples include nodal control volumes and surface areas in the finite volume case \cite{ref:NORD03}, Gauss-Lobatto quadrature weights for collocated spectral element operators \cite{carp96,gassner13}, as well as Gregory quadratures for high order finite difference operators \cite{ZING5,Linders18_acc}.  
For the case of structured methods on curvilinear domains, 
 extending the SBP operator formalism to the physical (as opposed to the reference) coordinate system has been the target of two recent papers \cite{Alund20_JCP, Lundquist18}.  In addition to this, the general SBP property (\ref{eq:SBP_partial0}) can be extended to a combination of multi-block, multi-element and hybrid meshes by encapsulating the numerical interface treatments into the global operator definition itself \cite{Lundquist22}.

The formal matrix structure required for the SBP relation (\ref{eq:SBP_partial0}) can be expressed as
\begin{equation}\label{eq:SBP_partial}
\mathcal{D}_{x_i}=\nwu{\mathcal{P}}^{-1}\mathcal{Q}_{x_i}, \quad  \mathcal{Q}_{x_i} +\mathcal{Q}_{x_i}^T= E^T \nwu{P} \nwu{N}_{i}E.
\end{equation}
In other words, each partial derivative operator $\mathcal{D}_{x_i}$ has an underlying almost skew-symmetric structure.  As already noted for $\mathcal{P}$, the specific values in $Q_{x_i}$ are determined by the particular numerical scheme.
\begin{remark}
\textcolor{red}{
The following semi-discrete analysis requires a diagonal matrix $\nwu{\mathcal{P}}$ in (\ref{eq:SBP_partial}). Therefore, it can unfortunately not be directly extended to standard finite element methods unless collocation/mass lumping is used.}
\end{remark}

\subsection{Partial time derivative approximation}

On a moving mesh, the dependent variables must of course follow the motion of the mesh nodes. Therefore,  any time derivative in the semi-discrete setting naturally corresponds to a material time derivative in the continuous setting. However, since the PDE problem (\ref{eq:cont1}) is formulated in terms of the partial derivative $u_t$ (corresponding to the original, physical coordinate system), we need to relate this term to the time derivative $d\boldsymbol{U}/dt$ of the numerical solution vector $\boldsymbol{U}(t)$.  \textcolor{red}{For this purpose, } we discretize the split form of the material derivative definition (\ref{eq:material_cont_split}) using a set of SBP operators (\ref{eq:SBP_partial}) in space.  \textcolor{red}{This leads to the following definition.}
\begin{definition}
\textcolor{red}{Consider a vector function $\boldsymbol{\Phi}(t)$ approximating $\phi\in\mathcal{U}$ on the moving mesh nodes,  and let $d\boldsymbol{\boldsymbol{\Phi}}/dt$ denote the time derivative of $\boldsymbol{\Phi}$. We then define the discrete partial time derivative $\boldsymbol{\Phi}_t$ approximating $\phi_t$ as}
\begin{equation}
\label{eq:material_disc_split}
 \textcolor{red}{\boldsymbol{\Phi}_t  =: \frac{d\boldsymbol{\Phi}}{dt}-}
\big(\mathcal{D}_m   - \frac{1}{2}\nwu{\boldsymbol{\nabla}\cdot \dot{\boldsymbol{X}}}\big) \boldsymbol{\Phi},
\end{equation}
where we employ a set of SBP (\ref{eq:SBP_partial}) partial derivative operators in order to define
\begin{equation}
\label{eq:Q_disc}
\mathcal{D}_m = \frac{1}{2}\Big(\nwu{\dot{\boldsymbol{X}}}_1\mathcal{D}_{x_1}+\ldots +\nwu{\dot{\boldsymbol{X}}}_p\mathcal{D}_{x_p} +\mathcal{D}_{x_1}\nwu{\dot{\boldsymbol{X}}_1} + \ldots  +\mathcal{D}_{x_p}\nwu{\dot{\boldsymbol{X}}}_p \Big),
\end{equation}
approximating the continuous operator $D_m$ in (\ref{eq:Q_cont}),
where the vectors $\dot{\boldsymbol{X}}_i$ contain the nodal values of $\dot{x}_i$.
Finally, the composite symbol $\boldsymbol{\nabla}\cdot\dot{\boldsymbol{X}} $ in (\ref{eq:material_disc_split}) denotes a vector containing either exact or approximated values of the continuous quantity $\nabla\cdot\dot{x}$.  Recall that the underscore notation in (\ref{eq:material_disc_split}) and (\ref{eq:Q_disc}) modifies the interpretation of these symbols from vectors to the corresponding diagonal matrices.  
\end{definition}

Since \textcolor{red}{$\nwu{\mathcal{P}}$}, $\dot{\boldsymbol{X}}_i$ as well as $E^T \nwu{P} \nwu{N}_{i}E$ are diagonal matrices (and thus commute), multiplying (\ref{eq:Q_disc}) from the left with \textcolor{red}{$\nwu{\mathcal{P}}$} and adding the transpose yields, using (\ref{eq:SBP_partial}),
\begin{equation*}
\begin{aligned}
\textcolor{red}{\nwu{\mathcal{P}}}\mathcal{D}_m+ \mathcal{D}_m^T\textcolor{red}{\nwu{\mathcal{P}}}= \frac{1}{2}\Big(
&
\nwu{\dot{\boldsymbol{X}}}_1E^T \nwu{P} \nwu{N}_{1}E +\nwu{\dot{\boldsymbol{X}}}_E^T \nwu{P} \nwu{N}_{2}E + \ldots +\nwu{\dot{\boldsymbol{X}}}_pE^T \nwu{P} \textcolor{red}{\nwu{N}_{p}}E 
+ \\ &  
 E^T \nwu{P} \nwu{N}_{1}E \nwu{\dot{\boldsymbol{X}}}_1+
 E^T \nwu{P} \nwu{N}_{2}E\nwu{\dot{\boldsymbol{X}}}_2 +\ldots +
 E^T \nwu{P} \textcolor{red}{\nwu{N}_{p}}E \nwu{\dot{\boldsymbol{X}}}_p
\Big) 
\\ = 
\Big(
&
\nwu{\dot{\boldsymbol{X}}}_1E^T \nwu{P} \nwu{N}_{1}E +\nwu{\dot{\boldsymbol{X}}}_2E^T \nwu{P} \nwu{N}_{2}E + \ldots +\nwu{\dot{\boldsymbol{X}}}_pE^T \nwu{P} \textcolor{red}{\nwu{N}_{p}}E\Big).
\end{aligned}
\end{equation*}
Note that each term $\nwu{\dot{\boldsymbol{X}}}_iE^T \nwu{P} \nwu{N}_{i}E$ above is only non-zero in diagonal positions corresponding to the surface nodes.  Hence,  multiplying from the left with $E^TE$ (i.e. a diagonal matrix with the constant value $1$ at the surface nodes) does not alter the value of this expression.  We can therefore write
\begin{equation*}
\nwu{\dot{\boldsymbol{X}}}_iE^T \nwu{P} \nwu{N}_{i}E = E^T\big(E\nwu{\dot{\boldsymbol{X}}}_iE^T \nwu{P} \nwu{N}_{i}\big)E  = E^T\nwu{\dot{X}}_i\nwu{P}\nwu{N}_i E =E^T\nwu{P}\nwu{N}_i\nwu{\dot{X}}_i E,
\end{equation*}
where $\nwu{\dot{X}}_i=E\nwu{\dot{\boldsymbol{X}}}_iE^T$ is the restriction of $\nwu{\dot{\boldsymbol{X}}}_i$ to the surface nodes.
Hence, the operator $\mathcal{D}_m$ in (\ref{eq:Q_disc}) satisfies the combined SBP property
\begin{equation}\label{eq:SBP_total_Q}
\mathcal{D}_m = \nwu{\mathcal{P}}^{-1}\mathcal{Q}_m, \quad \mathcal{Q}_m+\mathcal{Q}_m^T=E^T \nwu{P}
 \big(\nwu{N^T\dot{X}}\big)  E,
\end{equation}
where 
\begin{equation*}
\nwu{N^T\dot{X}} = :\nwu{N_{1}} \nwu{\dot{X}_1}+\nwu{N_{2}} \nwu{\dot{X}_2}+ \ldots +\nwu{N_{p}} \nwu{\dot{X}_p},
\end{equation*}
Note that (\ref{eq:SBP_total_Q}) leads to the inner product formulation of the same property
\begin{equation*}
\big( \boldsymbol{\Phi},\mathcal{D}_m\boldsymbol{\Psi}\big)_\mathcal{P} = -\big(  
\mathcal{D}_m\boldsymbol{\Phi},\boldsymbol{\Psi}  \big)_\mathcal{P}+
 \big(\Phi,  \big(\nwu{N^T\dot{X}}\big) \Psi \big)_P,
\end{equation*}
i.e. an exact analogue of (\ref{eq:IBP_material}).

\subsection{The semi-discrete scheme}
In order to impose the boundary condition in (\ref{eq:cont1}) in a transparent way, 
we take the lifting operator formulation (\ref{eq:cont1_lifting}) 
as basis for discretization.  However, we stress that this choice should not be viewed as a necessary feature in the current semi-discrete framework -- other means of imposing boundary conditions may also be considered.

We thus consider semi-discrete approximations to (\ref{eq:cont1}) satisfying the general simultaneous-approximation-term (SAT) \cite{ref:CARP94} format
\begin{equation}\label{eq:semi_disc1} 
\begin{aligned}
\boldsymbol{U}_t=\ &\mathcal{D}(\boldsymbol{U},t)\boldsymbol{U} +\mathcal{L}_{\boldsymbol{\delta}} \big[\mathcal{B}(\boldsymbol{U},t)\boldsymbol{U}-G(t)\big]+\boldsymbol{F}(t), && t\geq t_0 \\
\boldsymbol{U}(t_0) = \ & \boldsymbol{U}_0, 
\end{aligned}
\end{equation}
corresponding to the continuous formulation (\ref{eq:cont1_lifting}),
where $\boldsymbol{U}_t$ is the partial time derivative approximation of $\boldsymbol{U}$ defined by (\ref{eq:material_disc_split}). \textcolor{red}{ As opposed to the continuous setting, we do not make any a priori assumptions regarding the function space of admissible numerical solutions in (\ref{eq:semi_disc1}).
Thus,  we adopt the perspective of a general difference method by allowing $\boldsymbol{U}$ to be any vector of size $n_V$. In (\ref{eq:semi_disc1}),} we have also defined the discrete lifting operator $\mathcal{L}_{\boldsymbol{\delta}}$ as
\begin{equation*}
\mathcal{L}_{\boldsymbol{\delta}}= \nwu{\mathcal{P}}^{-1}\boldsymbol{\delta}^T \nwu{P},
\end{equation*}
where $\textcolor{green}{\boldsymbol{\delta}}$ is a $n_S\times n_V$ matrix approximating \textcolor{green}{the operator} $\delta$. Note that this definition of a discrete lifting operator implies, for any two vectors $\boldsymbol{\Phi}$ and $\Psi$ of dimension $n_V$ and $n_S$ respectively,
\begin{equation}\label{eq:lifting_disc}
\big(\boldsymbol{\Phi},   \mathcal{L}_{\boldsymbol{\delta}} \Psi \big)_\mathcal{P} = \big(\textcolor{green}{\boldsymbol{\delta}}\boldsymbol{\Phi},  \Psi\big)_P,
\end{equation}
thus mimicking (\ref{eq:lifting_cont}). 
\begin{remark}
With regards to the \textcolor{green}{example} boundary conditions discussed in section \ref{sec:example},  the lifting operator based on $\textcolor{green}{\delta =  \frac{1}{2}(a-\dot{x}_e)-\epsilon\partial_x}$ (Dirichlet case) can be approximated using $\textcolor{green}{\boldsymbol{\delta} =\frac{1}{2}(a-\dot{x}_e)E -\epsilon E\mathcal{D}_x}$, while $\delta = 1$ corresponds to $\textcolor{green}{\boldsymbol{\delta}} = E$.
\end{remark}
For future reference, the complete system matrix in (\ref{eq:semi_disc1}) is given by
\begin{equation}\label{eq:matrix_disc0}
\mathcal{M} = \mathcal{D} + \mathcal{L}_{\boldsymbol{\delta}}\mathcal{B},
\end{equation}
combining the differential operator and the boundary operator into one.

\subsection{Reference quadrature}
Relating to the stationary reference domain $\hat{\Omega}$, we assume that the volume quadrature $\mathcal{P}$ can be factored as $\mathcal{P}(t) = \nwu{\mathcal{J}(t)}\hat{\mathcal{P}} $, where $\hat{\mathcal{P}}$ is vector with constant reference quadrature weights, and the discrete Jacobian $\mathcal{J}(t) > 0$ satisfies the initial value problem
\begin{equation}\label{eq:GCL_disc} 
\begin{aligned}
\frac{d {\mathcal{J}}}{dt} = \ & \big( \nwu{\boldsymbol{\nabla}\cdot \dot{\boldsymbol{X}}}\big)\mathcal{J}, && t\geq t_0 \\
\mathcal{J}(t_0) = \ & \mathcal{J}_0>0,
\end{aligned}
\end{equation}
corresponding to (\ref{eq:GCL_cont}) in the continuous case. We have here defined $\boldsymbol{\nabla}\cdot \dot{\boldsymbol{X}}$ to be the same exact or approximated quantity previously employed in (\ref{eq:material_disc_split}).

Note that if both $\mathcal{J}$ and $\boldsymbol{\nabla}\cdot \dot{\boldsymbol{X}}$ contain exact nodal values of $J$ and $\nabla\cdot \dot{x}$ respectively, then (\ref{eq:GCL_disc}) of course follows directly from (\ref{eq:GCL_cont}). If either of them contains approximated values however, a bit more care is needed. For an extended discussion of the important special case of curvilinear coordinates, we refer to Appendix \ref{sec:curvi}.

\begin{remark}
As in the continuous case we can always make constructive use of (\ref{eq:GCL_disc}) in order to define the function $\mathcal{J}(t)$, even if \textcolor{green}{the} time dependent coordinate mapping is \textcolor{green}{not} explicitly known. In particular, we \textcolor{green}{can} do this by specifying the initial vector $\mathcal{J}_0$ to have the constant value $1$ everywhere, corresponding to the choice $\hat{\Omega}=\Omega(t_0)$, see Remark \ref{remark:curvi}.
\end{remark}
With the semi-discrete framework thus set up to mimic the essential properties of IBP (\ref{eq:IBP_partial}) as well as (\ref{eq:GCL_cont}),
we can now step by step mimic the continuous analysis, starting from Lemma \ref{lemma:materialderiv}. 

\begin{lemma}\label{lemma:materialderiv_disc}
For any differentiable function $\boldsymbol{\Phi}(t)$ of dimension $n_V$, 
the time derivative of the corresponding variable $\hat{\boldsymbol{\Phi}}= \nwu{\sqrt{\mathcal{J}}}\boldsymbol{\Phi}$ satisfies
\begin{equation}\label{eq:material_disc_modif}
\frac{d\hat{\boldsymbol{\Phi}}}{dt} =
\nwu{\sqrt{\mathcal{J}}}\big(\boldsymbol{\Phi}_t+\mathcal{D}_m\boldsymbol{\Phi}
\big),
\end{equation}
where $\boldsymbol{\Phi}_t$ is defined \textcolor{red}{in (\ref{eq:Q_disc})}, and where $\mathcal{J}(t)>0$ satisfies the auxiliary set of initial value problems (\ref{eq:GCL_disc}).
\begin{proof}
The proof of this result is analogous to the proof of Lemma \ref{lemma:materialderiv}. Taking the time derivative of $\hat{\boldsymbol{\Phi}} = \nwu{\sqrt{\mathcal{J}}}\boldsymbol{\Phi}$, we first get
\begin{equation}\label{eq:step1}
\frac{d\hat{\boldsymbol{\Phi}} }{dt} =  \nwu{\sqrt{\mathcal{J}}}\frac{d\boldsymbol{\Phi} }{dt} +\nwu{\frac{d\sqrt{\mathcal{J}}}{dt} }\boldsymbol{\Phi} .
\end{equation}
Moreover, since $\mathcal{J}>0$ we can rewrite (\ref{eq:GCL_disc}) into
\begin{equation} \label{eq:GCL_disc_sqrt}
\frac{d \sqrt{\mathcal{J}}}{dt} =  \frac{1}{2} \big( \nwu{\boldsymbol{\nabla}\cdot \dot{\boldsymbol{X}}}\big)
\sqrt{\mathcal{J}}.
\end{equation}
After inserting (\ref{eq:material_disc_split}) and (\ref{eq:GCL_disc_sqrt}) into (\ref{eq:step1}), we get
\begin{equation*}
\frac{d\hat{\boldsymbol{\Phi}} }{dt} =
\nwu{\sqrt{\mathcal{J}}}\Big[ \boldsymbol{\Phi}_t +
\big(\mathcal{D}_m   - \frac{1}{2}\nwu{\boldsymbol{\nabla}\cdot \dot{\boldsymbol{X}}}\big) \boldsymbol{\Phi}\Big]+\frac{1}{2} \big(\nwu{\boldsymbol{\nabla}\cdot \dot{\boldsymbol{X}}}\big)
\nwu{\sqrt{\mathcal{J}}}\boldsymbol{\Phi},
\end{equation*} 
i.e. (\ref{eq:material_disc_modif}) follows.
\end{proof}
\end{lemma}
In terms of the semi-discrete scheme (\ref{eq:semi_disc1}),  using Lemma \ref{lemma:materialderiv_disc}
the variable substitution $\hat{\boldsymbol{U}} = \nwu{\sqrt{\mathcal{J}} }\boldsymbol{U}$ now leads to the equivalent problem
\begin{equation}\label{eq:disc_MOL}
\begin{aligned}
\frac{d\hat{\boldsymbol{U}}}{dt} = \ &  \hat{\mathcal{D}}(\hat{\boldsymbol{U}},  t)\hat{\boldsymbol{U}} +\nwu{\sqrt{\mathcal{J}}}\mathcal{L}_{\boldsymbol{\delta}} \Big[\hat{\mathcal{B}}(\hat{\boldsymbol{U}},t)\hat{\boldsymbol{U}}-G(t)\Big]+\nwu{\sqrt{\mathcal{J}}}\boldsymbol{F}(t), && t\geq t_0 \\
\hat{\boldsymbol{U}}(t_0) = \ & \nwu{\sqrt{\mathcal{J}}}\boldsymbol{U}_0, 
\end{aligned}
\end{equation}
where
\begin{equation*}
\hat{\mathcal{D}}= \nwu{\sqrt{\mathcal{J}}} \big(\mathcal{D} +\mathcal{D}_m\big) \nwu{\sqrt{\mathcal{J}}} ^{-1}, \quad
\hat{\mathcal{B}}=  \mathcal{B}\nwu{\sqrt{\mathcal{J}}}^{-1}.
\end{equation*}
Just as in the corresponding continuous formulation (\ref{eq:cont_MOL}),  there is no indefinite source term appearing in (\ref{eq:disc_MOL}).  

The complete system matrix in (\ref{eq:disc_MOL}) is given by
\begin{equation}\label{eq:matrix_disc}
 \hat{\mathcal{M}}=\hat{\mathcal{D}}+\nwu{\sqrt{\mathcal{J}}}\mathcal{L}_{\boldsymbol{\delta}}\hat{\mathcal{B}}.
\end{equation}
Before extending \textcolor{green}{Proposition \ref{thm:semi-b} to the semi-discrete case}, we will also need a discrete version of the Reynolds transport theorem.

\subsection{A semi-discrete Reynolds transport theorem}\label{sec:Reynolds_disc}
Similarly to (\ref{eq:norm_equi_cont}) in the continuous analysis, for any two vectors $\boldsymbol{\Phi}$ and $\boldsymbol{\Psi}$ the substitution $\hat{\boldsymbol{\Phi}}= \nwu{\sqrt{\mathcal{J}}}\boldsymbol{\Phi}$ and $\hat{\boldsymbol{\Psi}}= \nwu{\sqrt{\mathcal{J}}}\boldsymbol{\Psi}$  leads directly to the identity
 \begin{equation*}
 \label{eq:norm_equi_disc}
 \big(\boldsymbol{\Phi}, \boldsymbol{\Psi}\big)_\mathcal{P}=  \big(\boldsymbol{\Phi}, \nwu{\mathcal{J}}\boldsymbol{\Psi}\big)_{\hat{\mathcal{P}}}= \big(\hat{\boldsymbol{\Phi}} , \hat{\boldsymbol{\Psi}}\big)_{\hat{\mathcal{P}}}.
 \end{equation*}
Since the quadrature weights in $\hat{\mathcal{P}}$ are constant, we can write the time derivative of this inner product as
 \begin{equation*}
\frac{d}{dt} \big( \boldsymbol{\Phi}, \boldsymbol{\Psi} \big)_{\mathcal{P}} = \frac{d}{d t}\big(\hat{\boldsymbol{\Phi}}, \hat{\boldsymbol{\Psi}}\big)_{\hat{\mathcal{P}}}  = \Big( \frac{d\hat{\boldsymbol{\Phi}}}{d t}, \hat{\boldsymbol{\Psi}}\Big)_{\hat{\mathcal{P}}} + \Big(\hat{\boldsymbol{\Phi}} ,\frac{d\hat{\boldsymbol{\Psi}}}{d t}\Big)_{\hat{\mathcal{P}}} .
 \end{equation*}
 Using  Lemma \ref{lemma:materialderiv_disc}, we can further write the first term on the right hand side as
 \begin{equation*}
 \Big( \frac{d\hat{\boldsymbol{\Phi}}}{d t}, \hat{\boldsymbol{\Psi}}\Big)_{\hat{\mathcal{P}}} =  \Big( \nwu{\sqrt{\mathcal{J}}}\boldsymbol{\Phi}_t, \hat{\boldsymbol{\Psi}}\Big)_{\hat{\mathcal{P}}} +  \Big( \nwu{\sqrt{\mathcal{J}}}\mathcal{D}_m\boldsymbol{\Phi}, \hat{\boldsymbol{\Psi}}\Big)_{\hat{\mathcal{P}}} =  \Big( \boldsymbol{\Phi}_t, \boldsymbol{\Psi}\Big)_{\mathcal{P}}+  \Big( \mathcal{D}_m\boldsymbol{\Phi} ,\boldsymbol{\Psi}\Big)_{\mathcal{P}},
 \end{equation*}
and the second term can of course be written analogously.  Finally applying the SBP property (\ref{eq:SBP_total_Q}), we have proven the semi-discrete Reynolds transport theorem
\begin{equation}\label{eq:Reynolds2_disc}
\frac{d}{dt} \big( \boldsymbol{\Phi}, \boldsymbol{\Psi} \big)_{\mathcal{P}}= 
\Big( \boldsymbol{\Phi}_t ,\boldsymbol{\Psi} \Big)_{\mathcal{P}} +
\Big( \boldsymbol{\Phi},  \boldsymbol{\Psi}_t\Big)_{\mathcal{P}} +
\big( \boldsymbol{\Phi}, \big(\nwu{N^T\dot{X}} \big)\boldsymbol{\Psi} \big)_P,
\end{equation}
which exactly corresponds to (\ref{eq:Reynolds2_cont}) in the continuous case.

\subsection{The energy method}

The Reynolds transport theorem (\ref{eq:Reynolds2_disc}) with $\boldsymbol{\Phi}= \boldsymbol{\Psi} = \boldsymbol{U}$ directly yields, corresponding to the continuous property (\ref{eq:en_cont}),
\begin{equation}\label{eq:en_disc}
\frac{d}{dt}\|\boldsymbol{U}\|_{\mathcal{P}}^2 
=2 \big( \boldsymbol{U}, \boldsymbol{U}_t \big)_{\mathcal{P}} +
\big( U,\big(\nwu{N^T\dot{X}} \big)    U \big)_P,
\end{equation}
where $U = E\boldsymbol{U}$ is the solution vector restricted to the surface nodes.  
After inserting the semi-discrete equation (\ref{eq:semi_disc1}) into (\ref{eq:en_disc}) and using the lifting operator property (\ref{eq:lifting_disc}), we further get
\begin{equation*}
\frac{d}{dt}\|\boldsymbol{U}\|^2_\mathcal{P} =2 \big(\boldsymbol{U}, \mathcal{D}\boldsymbol{U}\big)_\mathcal{P}+2\big(  \textcolor{green}{\boldsymbol{\delta}}\boldsymbol{U}, \mathcal{B}\boldsymbol{U}-G\big)_P+2\big(\boldsymbol{U}, \boldsymbol{F})_\mathcal{P}+\big( U, \big(\nwu{N^T\dot{X}} \big)   U \big)_P.
\end{equation*}
Thus, in order to guarantee an energy estimate of the type (\ref{eq:en_disc_stat}) for $\boldsymbol{F} = G = 0$,  the condition
\begin{equation}\label{eq:enboun_disc2}
\boldsymbol{\mathcal{E}}(\mathcal{D},\mathcal{B},\boldsymbol{\Phi}) \leq \alpha \|\boldsymbol{\Phi}\|^2_\mathcal{P},
\end{equation}
must hold for all vectors $\boldsymbol{\Phi}$ of size $n_V$, where
\begin{equation}\label{eq:eps_disc}
\boldsymbol{\mathcal{E}}(\mathcal{D},\mathcal{B},\boldsymbol{\Phi}) =\big(\boldsymbol{\Phi}, \mathcal{D}\boldsymbol{\Phi}\big)_\mathcal{P}+\big(\textcolor{green}{\boldsymbol{\delta}} \boldsymbol{\Phi}, \mathcal{B}\boldsymbol{\Phi}\big)_P+\frac{1}{2} \big( \Phi, \big(\nwu{N^T\dot{X}}  \big)   \Phi \big)_P,
\end{equation}
and where $\Phi =E\boldsymbol{\Phi}$. Note that (\ref{eq:enboun_disc2}) is analogous to the weak condition (\ref{eq:enboun_cont2}), sufficient for the strong condition (\ref{eq:enboun_cont1}) in the continuous analysis.  

Note that we can also express (\ref{eq:enboun_disc2}) in matrix form as
\begin{equation}\label{eq:enboun_disc3}
\nwu{\mathcal{P}} \mathcal{M}+\mathcal{M}^T\nwu{\mathcal{P}}+E^T\nwu{P}\nwu{N^T\dot{X}}E -\alpha \nwu{\mathcal{P}}\leq 0,
\end{equation}
where $\mathcal{M}$ is the system matrix of the problem given in (\ref{eq:matrix_disc0}).
We will therefore in general assume that the scheme can be designed in such a way that (\ref{eq:enboun_disc3}) holds.
We are now ready to formulate the main result of the semi-discrete analysis, relating (\ref{eq:enboun_disc3}) to a semi-boundedness result of the form (\ref{eq:semi-b_disc1}).
\begin{proposition}\label{thm:semi-b_disc}
The condition (\ref{eq:enboun_disc3}) is equivalent to the system matrix $\hat{\mathcal{M}}$ in (\ref{eq:matrix_disc}) satisfying
\begin{equation}\label{eq:semi-b_disc}
\nwu{\hat{\mathcal{P}}} \hat{\mathcal{M}} + \hat{\mathcal{M}}^T\nwu{\hat{\mathcal{P}}}  -\alpha \nwu{\hat{\mathcal{P}}}  \leq 0,
%\nwu{\hat{\mathcal{P}}} \big(\hat{\mathcal{D}}+\nwu{\sqrt{\mathcal{J}}}\mathcal{L}_{\boldsymbol{\delta}}\hat{\mathcal{B}}\big)+\big(\hat{\mathcal{D}}+\nwu{\sqrt{\mathcal{J}}}\mathcal{L}_{\boldsymbol{\delta}}\hat{\mathcal{B}}\big)^T\nwu{\hat{\mathcal{P}}}  -\alpha \hat{\mathcal{P}} \leq 0,
\end{equation}
thus guaranteeing that the real part of the matrix spectrum is non-positive if $\alpha \leq 0$.
\begin{proof}
First we can write, using the definition of $\hat{\mathcal{D}}$ in (\ref{eq:disc_MOL}),
\begin{equation*}
\nwu{\hat{\mathcal{P}}} \hat{\mathcal{M}} = \nwu{\hat{\mathcal{P}}}\nwu{\sqrt{\mathcal{J}}}\big(\mathcal{D} +\mathcal{L}_{\boldsymbol{\delta}}\mathcal{B}+\mathcal{D}_m\big)\nwu{\sqrt{\mathcal{J}}}^{-1}  =\nwu{\sqrt{\mathcal{J}}}^{-1} \nwu{\mathcal{P}}\big(\mathcal{M} +\mathcal{D}_m \big)\nwu{\sqrt{\mathcal{J}}}^{-1} .
\end{equation*}
By adding the transpose and using SBP (\ref{eq:SBP_total_Q}),  we now get
\begin{equation*}
\nwu{\hat{\mathcal{P}}} \hat{\mathcal{M}} + \hat{\mathcal{M}}^T\nwu{\hat{\mathcal{P}}}  = \nwu{\sqrt{\mathcal{J}}}^{-1} \Big[\nwu{\mathcal{P}} \mathcal{M} +\mathcal{M} ^T\nwu{\mathcal{P}}+E^T\nwu{P}\nwu{N^T\dot{X}}E\Big
]\nwu{\sqrt{\mathcal{J}}}^{-1} .
\end{equation*}
The whole left hand side of (\ref{eq:semi-b_disc}) can thus be written as, since $ \mathcal{P}=\nwu{\mathcal{J}}\hat{\mathcal{P}}$,
\begin{equation*}
\nwu{\hat{\mathcal{P}}} \hat{\mathcal{M}} + \hat{\mathcal{M}}^T\nwu{\hat{\mathcal{P}}}  -\alpha \nwu{\hat{\mathcal{P}}}  = \nwu{\sqrt{\mathcal{J}}}^{-1} \Big[\nwu{\mathcal{P}} \mathcal{M} +\mathcal{M} ^T\nwu{\mathcal{P}}+E^T\nwu{P}\nwu{N^T\dot{X}}E-\alpha \nwu{\mathcal{P}} \Big
]\nwu{\sqrt{\mathcal{J}}}^{-1} .
\end{equation*}
Finally, since $\nwu{\sqrt{\mathcal{J}}}^{-1}$ is a diagonal matrix containing strictly positive values on the diagonal, this shows that the two properties (\ref{eq:enboun_disc3}) and (\ref{eq:semi-b_disc}) are indeed equivalent.
\end{proof}
\end{proposition}
Since standard time integration methods require non-positive real parts of the semi-discrete matrix spectrum for stability,  the result of Proposition \ref{thm:semi-b_disc} is of fundamental importance for the method-of-lines approach to solving general PDE problems on a mesh subject to arbitrary motion;
if an energy estimate of the form (\ref{eq:enboun_disc3}) can be shown to hold for $\alpha \leq 0$, then (\ref{eq:disc_MOL}) can be marched forward in time in a stable way using standard (e.g. $A-$stable implicit or conditionally stable explicit) methods.

\section{Free-stream preservation}\label{sec:freestream}
In this section we investigate the additional often desired (but not essential, according to Proposition \ref{thm:semi-b_disc}) property of free-stream preservation. As we shall find, this property can be guaranteed by construction for both semi-discrete and fully discrete solutions using the ALE framework developed above.

\subsection{The semi-discrete case}
First, we require that the spatial approximation terms in (\ref{eq:semi_disc1}) are consistent with a constant solution $u_\infty$ in both space and time with the right choice of initial, boundary and forcing data.  We thus make the assumption \textcolor{red}{that for any real number $u_\infty$ we have}
\begin{equation}\label{eq:FSP_scheme}
\mathcal{D}(\boldsymbol{U}_\infty, t)\boldsymbol{U}_\infty = 0, \quad  \mathcal{B}(\boldsymbol{U}_\infty, t)\boldsymbol{U}_\infty = G, \quad \boldsymbol{F}(t)= 0,
\end{equation}
where $\boldsymbol{U}_\infty= u_\infty \mathbb{1}$, and $\mathbb{1}$ is a vector with the constant value $1$ in each position.  \textcolor{blue}{In the case of a constant coefficient problem,  the first condition in (\ref{eq:FSP_scheme}) follows directly from the use of consistent operators in space, see (\ref{eq:freestream_space}) below. }
Note that this leads to $\boldsymbol{U}_t= 0$ in (\ref{eq:semi_disc1}) after  inserting $\boldsymbol{U}=\boldsymbol{U}_\infty$ on the right hand side.  

Next,  recall that in the semi-discrete framework of section \ref{sec:semidisc}, the same vector $\boldsymbol{\nabla}\cdot \dot{\boldsymbol{X}}$ is used both in order to define $\boldsymbol{U}_t$ according to (\ref{eq:material_disc_split}) as well as in the auxiliary Jacobian equation (\ref{eq:GCL_disc}). For stability according to Proposition \ref{thm:semi-b_disc}, the choice of exact or approximated values of $\nabla\cdot\dot{x}$ was not essential. For the purpose of free-stream preservation however,  it is.  In particular, we consider the approximation obtained by applying the same set of SBP operators as previously for $\mathcal{D}_m$ in (\ref{eq:Q_disc}), i.e. we define
\begin{equation}\label{eq:DGCL}
\boldsymbol{\nabla}\cdot \dot{\boldsymbol{X}} = 
\mathcal{D}_{x_1}\dot{\boldsymbol{X}}_1+
\mathcal{D}_{x_2}\dot{\boldsymbol{X}}_2 +\ldots + \mathcal{D}_{x_p}\dot{\boldsymbol{X}}_p.
\end{equation}
Finally, we consider using consistent partial derivative approximations above, i.e. we assume that
\begin{equation}\label{eq:freestream_space}
\mathcal{D}_{x_i} \mathbb{1} = 0, \quad i=1,2,\ldots,p.
\end{equation}
\textcolor{blue}{
As shown in \ref{sec:curvi_FSP}, such consistency conditions can be achieved by construction for high order SBP operators using curvilinear coordinates.
However, as discussed in \cite{Lundquist20_JCP},  they are non-trivial (but still possible) to achieve in the presence of curved, non-collocated interfaces between computational blocks or elements.}
We can now prove
\begin{proposition}\label{thm:FSP_semidisc}
Consider using the same set of partial derivative operators $\mathcal{D}_{x_i}$  in (\ref{eq:Q_disc}) as well as in (\ref{eq:DGCL}),  and assume moreover that the consistency conditions (\ref{eq:FSP_scheme}) and (\ref{eq:freestream_space}) hold.
Then the semi-discrete scheme (\ref{eq:semi_disc1}), and thus equivalently (\ref{eq:disc_MOL}), is free-stream preserving.
\begin{proof}
Let the solution to (\ref{eq:semi_disc1}) at some point in time $t_1$ be given by the constant vector $\boldsymbol{U}(t_1)=\mathbb{1}$.  If both (\ref{eq:DGCL}) and (\ref{eq:freestream_space}) hold, then from the definition of $\mathcal{D}_m$ in (\ref{eq:Q_disc}) we have
\begin{equation*}
\mathcal{D}_m\mathbb{1} =  \frac{1}{2}\big(
\mathcal{D}_{x_1}\dot{\boldsymbol{X}}_1 +
\mathcal{D}_{x_2}\dot{\boldsymbol{X}}_2 +
\ldots
\big) = \frac{1}{2}
\boldsymbol{\nabla}\cdot\dot{\boldsymbol{X}}.
\end{equation*}
The time derivative of $\boldsymbol{U}$ at $t=t_1$ thus becomes, by inserting $\boldsymbol{U}=\boldsymbol{U}_\infty$ into the right hand side of (\ref{eq:material_disc_split}), 
\begin{equation*}
\frac{d\boldsymbol{U}}{dt}(t_1) = \boldsymbol{U}_t(t_1)+  \big(\mathcal{D}_m   - \frac{1}{2}\nwu{\boldsymbol{\nabla}\cdot \dot{\boldsymbol{X}}}\big)\boldsymbol{U}_\infty =
\boldsymbol{U}_t(t_1) +  \big(\frac{1}{2}
\boldsymbol{\nabla}\cdot\dot{\boldsymbol{X}} - \frac{1}{2}
\boldsymbol{\nabla}\cdot\dot{\boldsymbol{X}} \big)u_\infty= \boldsymbol{U}_t(t_1) .
\end{equation*}
As already discussed above,  $\boldsymbol{U}_t(t_1) = 0$ now follows from inserting (\ref{eq:FSP_scheme}) into the scheme (\ref{eq:semi_disc1}) with $\boldsymbol{U}=\boldsymbol{U}_\infty$.
\end{proof}
\end{proposition}

\subsection{Fully discrete approximations}
The two schemes (\ref{eq:semi_disc1}) and (\ref{eq:disc_MOL}) are equivalent as long as the time variable is kept continuous.  When marching the system forward in time however, applying the same time integration method with the same time step size will in general produce two different numerical solutions.  Since Proposition \ref{thm:semi-b_disc} can be used to infer the time-stability of (\ref{eq:disc_MOL}) as long as $\alpha\leq 0$,  this formulation is most suitable to integrate in time. In order to obtain free-stream preservation in a fully discrete sense, we thus need to verify that the result of Proposition \ref{thm:FSP_semidisc} extends to solutions obtained by marching (\ref{eq:disc_MOL}) forward in time using arbitrary time step sizes.  
To this end, an approach based on updating the scheme as well as the metric Jacobian \textcolor{green}{determinant} using the same time integration scheme has been considered in some previous works,  see e.g.  \textcolor{blue}{\cite{Persson09, Minoli11, Badia06,Bonito13a,Bonito13b}}. We will consider an appropriate modification to this approach below based on the squared root relation (\ref{eq:disc_MOL}).

In terms of the new independent variable $\hat{\boldsymbol{U}}=\nwu{\sqrt{\mathcal{J}}}\boldsymbol{U}$,  a constant solution $ \boldsymbol{U}=\boldsymbol{U}_\infty$ is of course equivalent to $ \hat{\boldsymbol{U}}=u_\infty\sqrt{\mathcal{J}}$.  In order to retain free-stream preservation after marching (\ref{eq:disc_MOL}) forward in time, we need $\hat{\boldsymbol{U}}=u_\infty\sqrt{\mathcal{J}}$ to hold for fully discrete solutions as well. 
Thus, consider (\ref{eq:GCL_disc_sqrt}) together with (\ref{eq:disc_MOL}) as the semi-coupled system
\begin{equation}\label{eq:system}
\begin{aligned}
\frac{d \sqrt{\mathcal{J}}}{dt} =  \ & \frac{1}{2} \big( \nwu{\boldsymbol{\nabla}\cdot \dot{\boldsymbol{X}}}\big)
\sqrt{\mathcal{J}} \\
\frac{d\hat{\boldsymbol{U}}}{dt} = \ &  RHS\big(\sqrt{\mathcal{J}}, \boldsymbol{U}, t\big),
\end{aligned}
\end{equation}
where we have introduced the shorthand notation
\begin{equation*}
RHS\big(\sqrt{\mathcal{J}}, \boldsymbol{U}, t\big)= \nwu{\sqrt{\mathcal{J}}} \Big[\mathcal{D}_m\boldsymbol{U}+ \mathcal{D}\big(\boldsymbol{U},t\big)\boldsymbol{U} + \mathcal{L} \big[\mathcal{B}(\boldsymbol{U},t)\boldsymbol{U}-G(t)\big]+\boldsymbol{F}(t)\Big] .
\end{equation*}
Notice that
\begin{equation}\label{eq:FSP_sqrt2}
RHS\big(\sqrt{\mathcal{J}}, \textcolor{green}{\boldsymbol{U}_\infty}, t\big)= \frac{1}{2} \big( \nwu{\boldsymbol{\nabla}\cdot \dot{\boldsymbol{X}}}\big)
u_\infty\sqrt{\mathcal{J}} ,
\end{equation}
following the same consistency assumptions as in Proposition \ref{thm:FSP_semidisc}.
Thus, by starting from a constant solution $\boldsymbol{U}(t_1) = \boldsymbol{U}_\infty$, the right hand side of the second equation in (\ref{eq:system}) reduces to exactly the right hand side of the first equation, only scaled by $u_\infty$. As long as the same time integration scheme is applied to both equations,  we thus automatically obtain $\hat{\boldsymbol{U}}=u_\infty\sqrt{\mathcal{J}}$ for all additional discrete time steps.  

To illustrate this \textcolor{green}{result} in more explicit terms, we give the full proof of free-stream preservation below for the special case of applying a general $s-$stage explicit Runge-Kutta time marching scheme.  

\begin{proposition}\label{thm:RK}
Consider marching the system (\ref{eq:system}) forward in time using a general $s-$stage explicit Runge-Kutta time marching scheme, i.e. we apply
\begin{equation}\label{eq:RK_scheme}
\begin{aligned}
\sqrt{{\mathcal{J}}^{n+1}} = & \ \sqrt{{\mathcal{J}}^{n}} + \Delta t \sum_{k=1}^s 
 \frac{b_k}{2} \big( \nwu{\boldsymbol{\nabla}\cdot \dot{\boldsymbol{X}}}\big(t_n+c_k\Delta t\big)\big)
\sqrt{\bar{\mathcal{J}}^k} \\
\nwu{\sqrt{{\mathcal{J}}^{n+1}}}\boldsymbol{U}^{n+1} = & \ \nwu{\sqrt{{\mathcal{J}}^n}}\boldsymbol{U}^n + \Delta t \sum_{k=1}^s 
b_k RHS\big(\sqrt{\bar{\mathcal{J}}^k},\bar{\boldsymbol{U}}^k,t_n+c_k\Delta t\big),
\end{aligned}
\end{equation}
where the stage values $\bar{\mathcal{J}}^k$ and $\bar{\boldsymbol{U}}^k$ are defined by
\begin{equation}\label{eq:RK_stage}
\begin{aligned}
\sqrt{\bar{\mathcal{J}}^{k}} = & \ \sqrt{{\mathcal{J}}^{n}} +
 \Delta t \sum_{\nu=1}^{k-1} 
 \frac{a_{k\nu}}{2}
 \big( \nwu{\boldsymbol{\nabla}\cdot \dot{\boldsymbol{X}}}\big(t_n+c_\nu\Delta t\big)\big)\sqrt{\bar{\mathcal{J}}^\nu}  \\
\nwu{\sqrt{\bar{\mathcal{J}}^{k}}}\bar{\boldsymbol{U}}^{k} = & \ \nwu{\sqrt{{\mathcal{J}}^n}}\boldsymbol{U}^n+ \Delta t \sum_{\nu=1}^{k-1}
a_{k\nu} RHS\big(\sqrt{\bar{\mathcal{J}}^\nu},\bar{\boldsymbol{U}}^\nu,t_n+c_\nu \Delta t\big).
\end{aligned}
\end{equation}
If the premise of Proposition \ref{thm:FSP_semidisc} is satisfied, then this time marching scheme is free-stream preserving.
\begin{proof}
Starting from a constant solution $\boldsymbol{U}^n = \boldsymbol{U}_\infty=u_\infty \mathbb{1}$ we need to demonstrate that this automatically implies $\boldsymbol{U}^{n+1} =\boldsymbol{U}_\infty$, and we will do so by the method of induction over the stage index $k$. Thus, for $k=1$ in (\ref{eq:RK_stage}), $\boldsymbol{U}^n =\boldsymbol{U}_\infty$ yields
\begin{equation*}
\begin{aligned}
\sqrt{\bar{\mathcal{J}}^{1}} = \ & \sqrt{{\mathcal{J}}^{n}}, \\ 
\nwu{\sqrt{\bar{\mathcal{J}}^{1}}}\bar{\boldsymbol{U}}^{1} = & \ \nwu{\sqrt{{\mathcal{J}}^n}}\boldsymbol{U}^n =u_\infty \sqrt{{\mathcal{J}}^{n}} = u_\infty\sqrt{\bar{\mathcal{J}}^{1}},
\end{aligned}
\end{equation*}
thus showing that $\bar{\boldsymbol{U}}^1 = \boldsymbol{U}_\infty$. 

Next,  we assume that $\bar{\boldsymbol{U}}^\nu =\boldsymbol{U}_\infty$ holds for $\nu = 1,\ldots , k -1 $, where $k$ is some number between $1$ and $s$.  The second equation in (\ref{eq:RK_stage}) then becomes, using (\ref{eq:FSP_sqrt2}), 
\begin{equation*}
\nwu{\sqrt{\bar{\mathcal{J}}^{k}}}\bar{\boldsymbol{U}}^{k} =  u_\infty\sqrt{{\mathcal{J}}^{n}} +
 \Delta t \sum_{\nu=1}^{k-1} 
 \frac{a_{k\nu}}{2}
 \big( \nwu{\boldsymbol{\nabla}\cdot \dot{\boldsymbol{X}}}\big(t_n+c_\nu\Delta t\big)\big)u_\infty\sqrt{\bar{\mathcal{J}}^\nu} =u_\infty \sqrt{\bar{\mathcal{J}}^{k}},
\end{equation*}
where in the last step we have used the first equation in (\ref{eq:RK_stage}).

We have thus shown that $\bar{\boldsymbol{U}}^{\nu} =\boldsymbol{U}_\infty$ for $\nu=1,...,k-1$ implies $\bar{\boldsymbol{U}}^k = \boldsymbol{U}_\infty$, and by induction it follows that $\bar{\boldsymbol{U}}^k =\boldsymbol{U}_\infty$ holds true for all stages $k=1,\ldots,s$. Inserting this result into the second line of the Runge-Kutta scheme (\ref{eq:RK_scheme}) now yields
\begin{equation*}
\nwu{\sqrt{{\mathcal{J}}^{n+1}}}\boldsymbol{U}^{n+1}  =u_\infty \sqrt{{\mathcal{J}}^{n}} + \Delta t \sum_{k=1}^s 
 \frac{b_k}{2} \big( \nwu{\boldsymbol{\nabla}\cdot \dot{\boldsymbol{X}}}\big(t_n+c_k\Delta t\big)\big)u_\infty
\sqrt{\bar{\mathcal{J}}^k}  = u_\infty\sqrt{{\mathcal{J}}^{n+1}},
\end{equation*}
i.e. we have proven that $\boldsymbol{U}^{n+1} =\boldsymbol{U}_\infty$.
\end{proof}
\end{proposition}

\section{\textcolor{red}{Numerical experiment}}
\label{sec:num}
\textcolor{red}{
In this section we validate the stability and accuracy of the proposed method-of-lines ALE technique  through a numerical experiment.
As a model, consider the scalar advection-diffusion equation in one dimension with a characteristic inflow condition on the left boundary $x_s(t)$, and a Dirichlet condition on the right boundary $x_e(t)$,
 \begin{equation}\label{eq:model_problem}
 \begin{aligned}
 u_t  = \ & \epsilon u_{xx} - u_x + F(x,t),  && x_s(t)< x<x_e(t),  &&t>0 \\
\epsilon u_x -(1-\dot{x}_s)u = \ & g_s(t), && x=x_s(t),&& t> 0 \\
u = \ & 0 && x=x_e(t) && t> 0 \\
u(x, 0) = \ &u_0(x),  &&  x_s< x<x_e.
 \end{aligned}
 \end{equation}
For simplicity, assume that $\dot{x}_s (t)\leq 1$ holds, so that the left boundary remains an inflow boundary at all times.  See also section \ref{sec:example} for an energy analysis of the two types of boundary conditions present in (\ref{eq:model_problem}).
Finally, we impose the time periodic manufactured solution $u= u_{\mathrm{MMS}}= (1 + A\mathrm{sin}(x-t))(1-e^{(x-x_e) /\epsilon})$, with $A=0.1$, by specifying the data
\begin{equation*}
\begin{aligned}
F = \ & A\epsilon \mathrm{sin}(x-t)+ \Big[\frac{\dot{x}_e}{\epsilon} (1+A\mathrm{sin}(x-t)) + 2A\mathrm{cos}(x-t)-A\epsilon\mathrm{sin}(x-t)\Big]e^{\frac{x-x_e}{\epsilon}} \\
% g_s = \ & 1+A\mathrm{sin}(x_s-t) - A\epsilon \mathrm{cos}(x_s-t)(1-e^{\frac{x_s-x_e}{\epsilon}})\\
 g_s = \ &  -(1+A\mathrm{sin}(x_s-t)) + \big[ A\epsilon \mathrm{cos}(x_s-t)+\dot{x}_s (1+A\mathrm{sin}(x_s-t))\big](1-e^{(x_s-x_e) /\epsilon})\\
 u_0 = \ & (1+A\mathrm{sin}(x)) (1-e^{(x-x_e(0)) /\epsilon}), 
\end{aligned}
\end{equation*}
}

\begin{figure}[htb!]
\begin{center}
\begin{tabular}{ll}
\includegraphics[scale=.53]{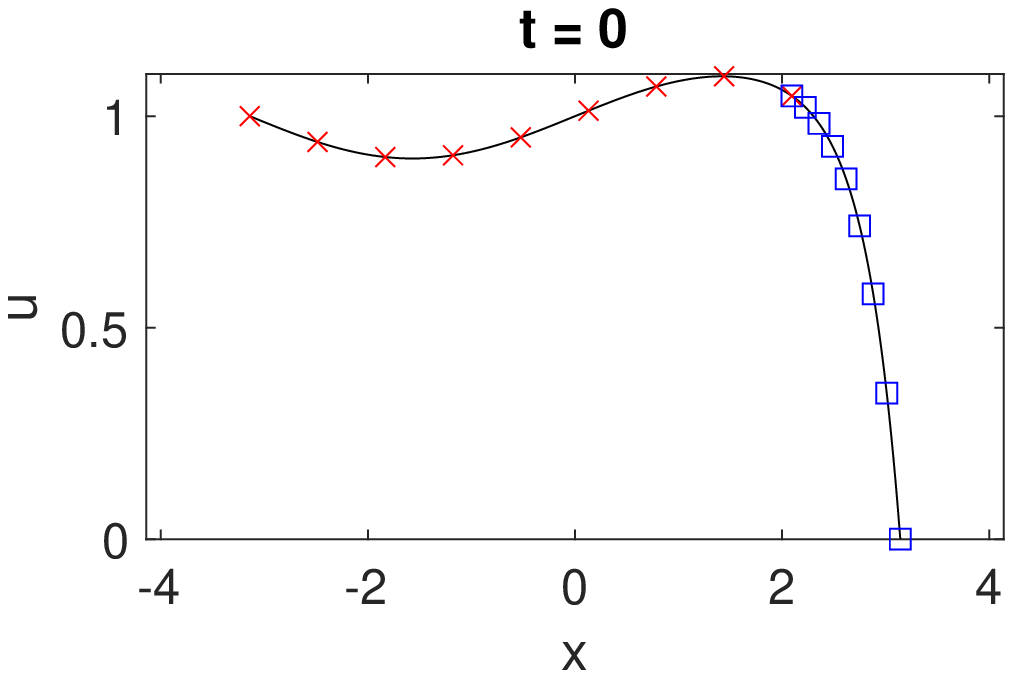} & \includegraphics[scale=.53]{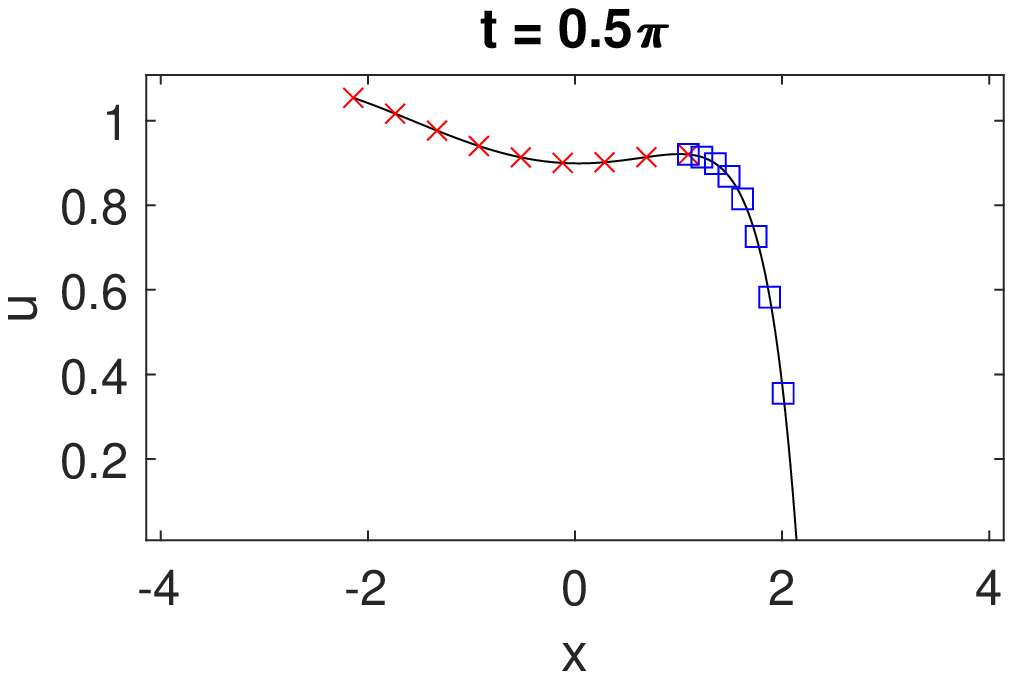} \\
\includegraphics[scale=.53]{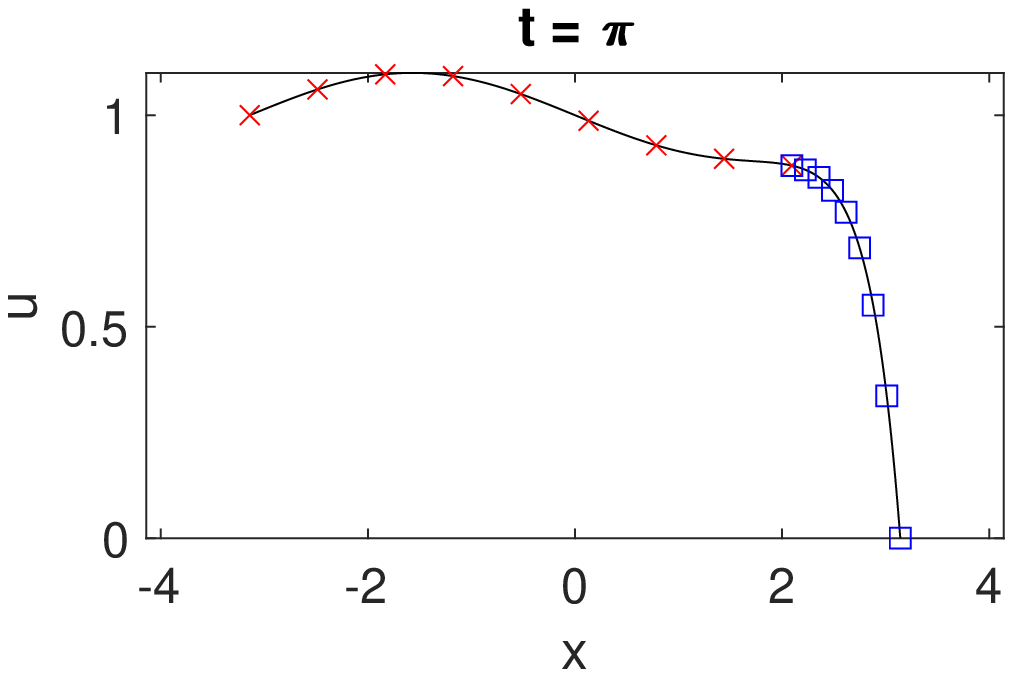} & \includegraphics[scale=.53]{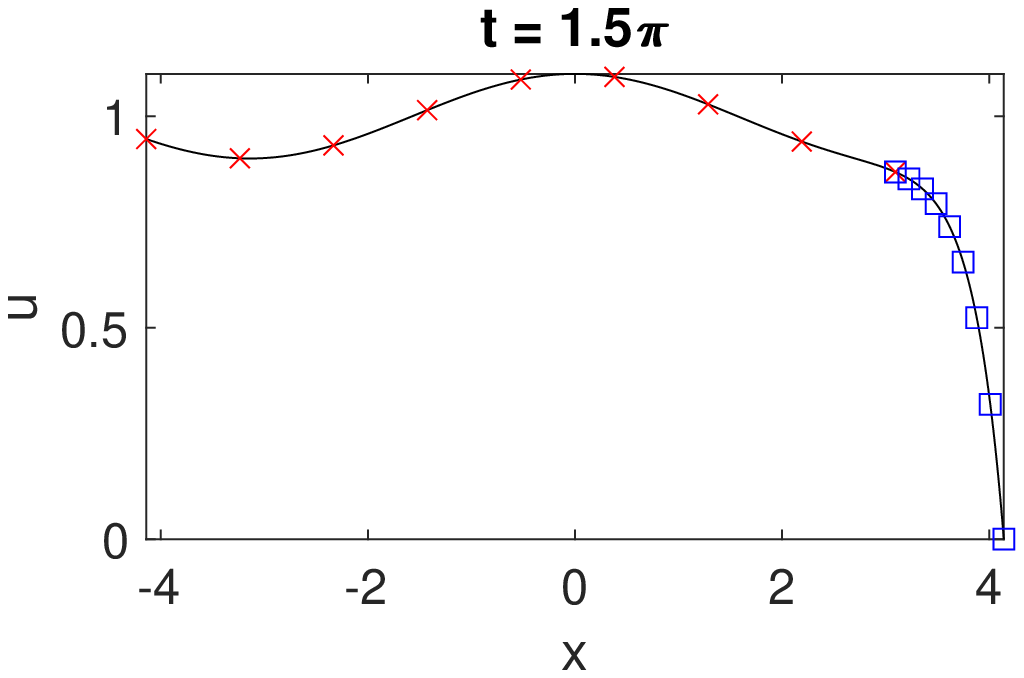}
\end{tabular}
\caption{Exact solution on the two subdomains during one periodic cycle. The ratio of grid spacings between the left and right domain is given by $5$ at $t=0$.}
\label{fig:exact} 
\end{center}
\end{figure}

\textcolor{red}{For the numerical experiment we specify $\epsilon = 0.1\pi$ and set the moving domain boundary as $x_s=-\pi+\mathrm{sin }(t)$, $x_e = \pi-\mathrm{sin}(t)$.  Furthermore, the domain is subdivided into two blocks given by $[x_s,x_e] = [x_s,x_m]\cup[x_m,x_e]$, where $x_m = x_e-\pi/3$. In other words, we insert a boundary layer grid of width $\pi/3$ in the rightmost part of the domain, undergoing rigid motion. Note that this combination of rigid and non-rigid motion in the two blocks leads to locally non-smooth mesh velocities at the block interface.  The ratio of grid spacings between the blocks is given by $5$ at $t=0$.  See Figure \ref{fig:exact} for an illustration of the moving block grid together with the exact solution to the model problem using $8$ grid spacings in the boundary layer.  Furthermore, we apply fourth order accurate finite difference stencils in the block interiors, amended with second order accurate SBP boundary closures, see e.g.  \cite{Gustafsson08} for the tabulated values.  Finally, we apply a SAT interface condition as in \cite{ref:NORD01,Lundquist22} in order to construct an encapsulated SBP operator on the full spatial domain. }

\textcolor{red}{On the discrete side, the SBP property (\ref{eq:SBP_partial}) in one space dimension simplifies to $\mathcal{D}_x=\mathcal{P}^{-1}\mathcal{Q}_x$, where
 \begin{equation}
\mathcal{Q}_x+\mathcal{Q}_x^T = -E_s^TE_s+E_e^TE_e \quad E_s = \begin{pmatrix}
1 & 0 & \ldots & 0
\end{pmatrix}, \quad E_e = \begin{pmatrix}
0 & \ldots & 0 & 1
\end{pmatrix}.
 \end{equation}
Based on the continuous analysis in (\ref{sec:example}), we specify the discrete operators in (\ref{eq:semi_disc1}) as 
\begin{equation*}
\mathcal{D} = \epsilon\mathcal{D}_x^2 - \mathcal{D}_x, \quad \mathcal{B} = \begin{pmatrix}
\epsilon E_s \mathcal{D}_x -(1-\dot{x}_s) E_s \\ E_e
\end{pmatrix}, \quad \boldsymbol{\delta} =  \begin{pmatrix}
E_s \\ -\epsilon E_e\mathcal{D}_x + \frac{1}{2}(1 - \dot{x}_e)E_e
\end{pmatrix},
\end{equation*}
where for simplicity we have applied the first derivative operator twice in order to approximate the second derivative.  In order to limit the resulting high wavenumber errors, we apply energy stable high order numerical filters \cite{Lundquist20_JSC} after each time step.
In (\ref{eq:eps_disc}), the above operators yield
\begin{equation*}
\big(\boldsymbol{\Phi}, \mathcal{D}\boldsymbol{\Phi}\big)_\mathcal{P} = -\epsilon\|\mathcal{D}_x\boldsymbol{\Phi}\|_\mathcal{P}^2 - \epsilon\Phi_s\Phi_{sx} + \frac{1}{2} \Phi_s^2+\epsilon \Phi_e\Phi_{ex} - \frac{1}{2} \Phi_e^2,
\end{equation*}
where $\Phi_s = E_s\boldsymbol{\Phi}$, $\Phi_e = E_e\boldsymbol{\Phi}$, $\Phi_{xs} = E_s\mathcal{D}_x\boldsymbol{\Phi}$, $\Phi_{xe} = E_e\mathcal{D}_x\boldsymbol{\Phi}$, as well as 
\begin{equation*}
\big(\boldsymbol{\delta} \boldsymbol{\Phi}, \mathcal{B}\boldsymbol{\Phi}\big)_P = \Phi_s \big[\epsilon \Phi_{xs} - (1-\dot{x}_s)\Phi_s\big] + \big[-\epsilon \Phi_{xe} + \frac{1}{2}(1-\dot{x}_e)\Phi_e\big] \Phi_e.
\end{equation*}
The full expression in (\ref{eq:eps_disc}) thus becomes
\begin{equation*}
\begin{aligned}
\boldsymbol{\mathcal{E}}(\mathcal{D},\mathcal{B},\boldsymbol{\Phi}) =  \ & -\epsilon\|\mathcal{D}_x\boldsymbol{\Phi}\|_\mathcal{P}^2 - \epsilon\Phi_s\Phi_{sx} - \frac{1}{2} (\dot{x}_s-1)\Phi_s^2+ \Phi_s \big[\epsilon \Phi_{xs} - (1-\dot{x}_s)\Phi_s\big] \\ \ & +\epsilon \Phi_e\Phi_{ex} + \frac{1}{2} (\dot{x}_e-1)\Phi_e^2+\big[-\epsilon \Phi_{xe}  +\frac{1}{2}(1-\dot{x}_e)\Phi_e\big] \Phi_e \\
= \ & -\epsilon\|\mathcal{D}_x\boldsymbol{\Phi}\|_\mathcal{P}^2 -\frac{1}{2}(1-\dot{x}_s)\Phi_s^2 \leq 0,
\end{aligned}
\end{equation*}
for all vectors $\boldsymbol{\Phi}$, showing that the semi-discrete scheme (\ref{eq:semi_disc1}) is energy stable.}

\begin{figure}[htb!]
\begin{center}
\includegraphics[scale=.53]{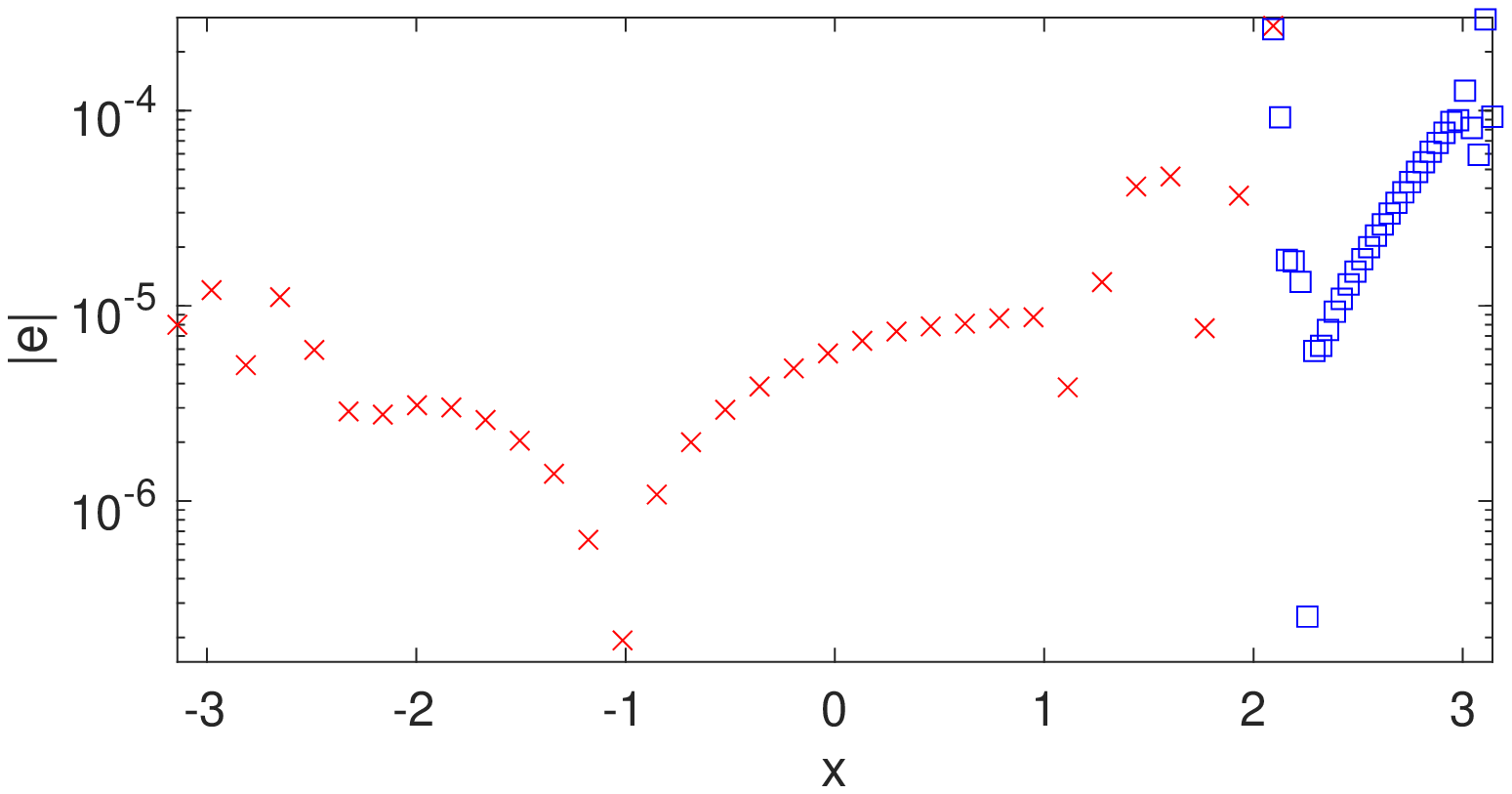} \\
\includegraphics[scale=.53]{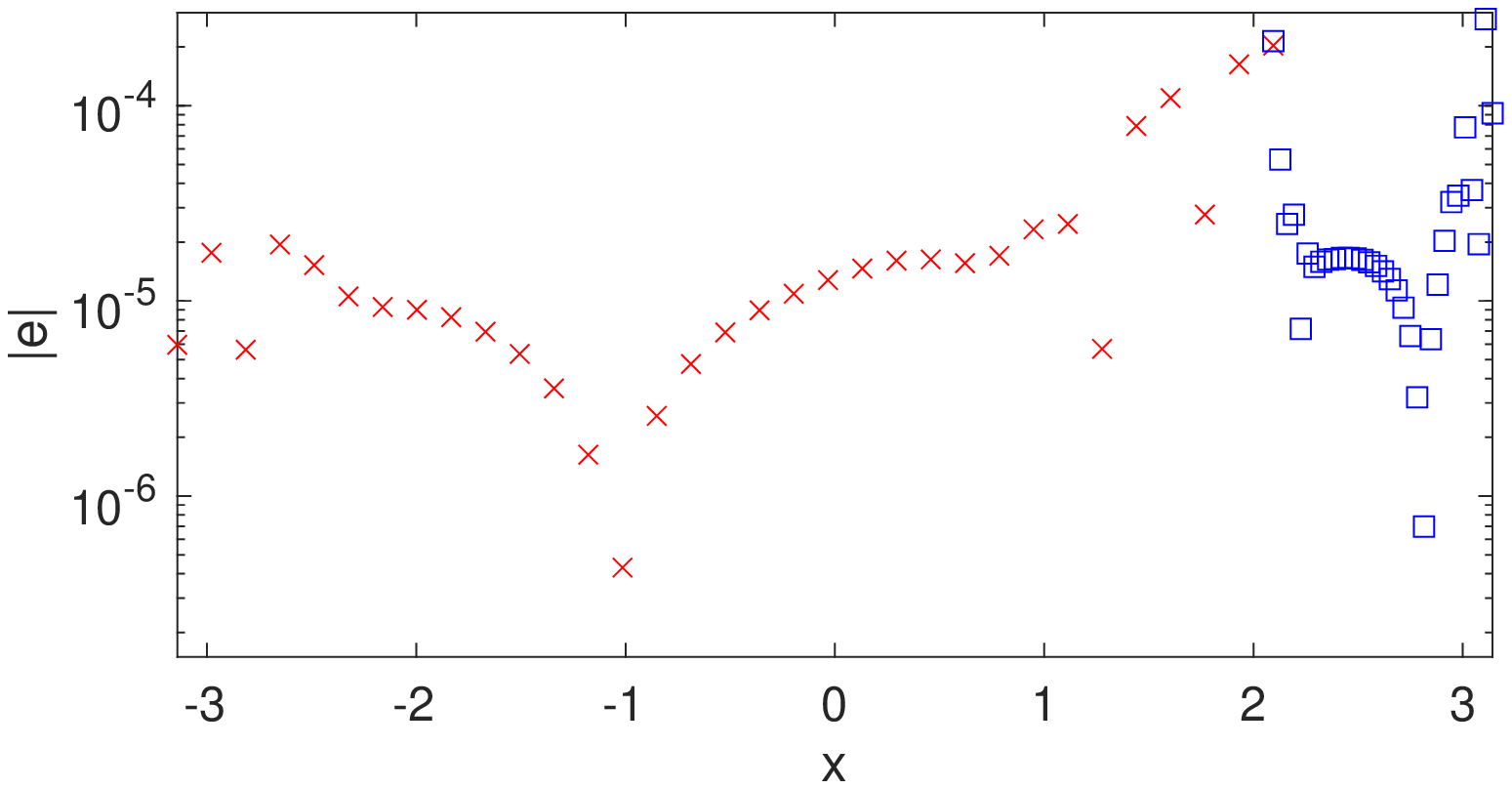}
\caption{Absolute error levels for the stationary (above) and moving (below) domain case at $t=2\pi$.}
\label{fig:errplot} 
\end{center}
\end{figure}

\textcolor{red}{
For time discretization we apply the classical fourth order explicit Runge-Kutta scheme with a time step size given by $\Delta t = 0.5 \Delta x_{\mathrm{min}}^2/\epsilon$.  To illustrate the fact that discrete geometric conservation is not necessary, we use the exact values of $\sqrt{J}$ in the scheme (\ref{eq:disc_MOL}) rather than solving the coupled system (\ref{eq:system}) in time. For comparison, we also consider the stationary domain case $x_s=-\pi$, $x_e=\pi$ in (\ref{eq:model_problem}) using the same block grid as in the moving domain case at $t=0$, as well as the same manufactured exact solution. 
The resulting absolute error levels at $t= 2\pi$ are plotted in Figure \ref{fig:errplot} using $N=32$ grid spacings in the boundary layer grid.
In Table \ref{table:conv} we also compare the convergence rates between the moving and stationary domain cases.  The convergence rate is close to three in both cases, with small differences in absolute error levels. The result of this experiment thus indicates that the proposed method-of-lines approach is comparable in terms of numerical efficiency between a moving domain and the equivalent stationary domain application.} 


\begin{table*}[htb!]\centering
\begin{tabular}{@{}ccccccccccc@{}}
& \multicolumn{2}{c}{stationary} & \phantom{abc}& \multicolumn{2}{c}{moving}\\
% &\phantom{abc} & \multicolumn{2}{c}{$s=3$}\\ 
 \cmidrule{2-3} \cmidrule{5-6} \cmidrule{8-9}
$N$ & $\mathrm{log}_{10}(\|e\|_\infty)$ & $p$ && $\mathrm{log}_{10}(\|e\|_\infty)$ & $p$
\\ \midrule
     8  & -1.8639  &    -      &&        -2.0192      &    -    &&         \\
    16 & -2.6741  &  2.69  && -2.7303  &    2.36    &&        \\
    32 & -3.5329  &  2.85  && -3.5558  & 2.74  &&  \\
    64 & -4.3916  &  2.85  && -4.4243   & 2.89  &&  \\
        128 & -5.2678  &  2.91  && -5.3045   & 2.92  &&  \\
 \bottomrule
\end{tabular}
\caption{Convergence in $L_\infty$ norm at $t=2\pi$ with $N$ the number of spacings in the boundary layer grid.}
\label{table:conv} 
\end{table*}

% 
% \begin{equation*}
%\begin{aligned}
%\boldsymbol{U}_t=\ & \epsilon \mathcal{D}_x^2\boldsymbol{U} +\mathcal{L}_{\boldsymbol{\delta}} \big[\mathcal{B}(\boldsymbol{U},t)\boldsymbol{U}-G(t)\big]+\boldsymbol{F}(t), && t\geq t_0 \\
%\boldsymbol{U}(t_0) = \ & \boldsymbol{U}_0, 
%\end{aligned}
%\end{equation*}

% \begin{equation*}
% \begin{aligned}
% u_t  = \ & \epsilon u_{xx} - au_x + F(x,t),  && 0 < x<1,  &&t>0 \\
%u(1,t) = \ & 0 && x =1, && t> 0 \\
%u(x, 0) = \ &u_0(x),  && x < x_e(t), 
% \end{aligned}
% \end{equation*}

\section{Conclusions}\label{sec:concl} 
We have developed a framework based on semi-bounded operators for analyzing and discretizing initial boundary value problems on domains and meshes subject to arbitrary motion. This framework can be used to both analyze the linear well-posedness of the continuous problem as well as to formulate energy stable semi-discretizations. In particular, we have demonstrated that any energy estimate derived with respect to the physical coordinate system is associated with a semi-bounded property on the stationary reference domain. Since, with the right set-up, the main continuous result could be proven based on integration-by-parts \textcolor{green}{only},  the semi-discrete framework could be developed at a high level of abstraction based on discrete summation-by-parts, and thus valid for many different types of numerical discretization methods.

As a result of this development, initial boundary value problems posed on moving domains can be discretized in space and time separately in a stable way, i.e. using a standard method-of-lines approach. This does not require that discrete geometric conservation (free-stream preservation) holds for the combined spatial and temporal discretization. As a final bonus we also demonstrated that free-stream preservation can be achieved if so desired, both on the semi-discrete level as well as for fully discrete solutions. However, we stress that this is not strictly necessary for either high order accuracy or stability. \textcolor{red}{The stability and efficiency of the new method-of-lines framework, even in the absence of discrete geometric conservation, was corroborated through a numerical experiment.}

\appendix

\section{Curvilinear coordinates}\label{sec:curvi}
Here we discuss in more detail the special case where the moving domain coordinates $(x_1, \ x_2,  \ \ldots, \ x_p)\in\Omega$ for each time $t$ are expressed by a smooth bijective mapping $x_i = x_i(t, \ \xi_1, \  \xi_2, \ \ldots, \ \xi_p )$ from the unit hypercube $\hat{\Omega}= [0,1]^p$. Provided that a discrete quadrature rule $\hat{\mathcal{P}}$ is available for $\hat{\Omega}$, we can then always define a corresponding quadrature for $\Omega$ as before using $\mathcal{P}(t) = \nwu{\mathcal{J}(t)}\hat{\mathcal{P}} $,  where $\mathcal{J}$ contains nodal values of the Jacobian determinant,  either exact or approximated. 

Note that any specific choice of $\nwu{\mathcal{J}(t)}$ simultaneously dictates the value of $\boldsymbol{\nabla}\cdot\dot{\boldsymbol{X}}
=  \nwu{\mathcal{J}}^{-1}\frac{d\mathcal{J}}{dt}$ by the relation (\ref{eq:GCL_disc}).
For example, in $2$D we have the exact expression
\begin{equation*}
J = 
\frac{\partial x_1}{\partial \xi_1}\frac{\partial x_2}{\partial \xi_2} - 
\frac{\partial x_1}{\partial \xi_2}\frac{\partial x_2}{\partial \xi_1} ,
\end{equation*}
and it is therefore natural to either insert these continuous values in $\mathcal{J}$ directly, or alternatively to apply discrete partial derivative operators on $\hat{\Omega}$ in order to define
\begin{equation}\label{eq:J_appr}
\mathcal{J} = 
\big(\nwu{\mathcal{D}_{\xi_1}\boldsymbol{X}_1}\big) 
\big(\mathcal{D}_{\xi_2}\boldsymbol{X}_2\big)  - 
\big(\nwu{\mathcal{D}_{\xi_2}\boldsymbol{X}_1}\big)  
\big(\mathcal{D}_{\xi_1}\boldsymbol{X}_2\big) .
\end{equation}
With this choice, the discrete variable $\boldsymbol{\nabla}\cdot\dot{\boldsymbol{X}}
=  \nwu{\mathcal{J}}^{-1}\frac{d\mathcal{J}}{dt}$ becomes
\begin{equation}\label{eq:J_FSP}
\begin{aligned}
\boldsymbol{\nabla}\cdot\dot{\boldsymbol{X}}
%= : \mathcal{J}^{-1}\frac{d\mathcal{J}}{dt}
=\nwu{\mathcal{J}}^{-1}
\Big[&
\big(\nwu{\mathcal{D}_{\xi_1}\dot{\boldsymbol{X}}_1}\big)
\big(\mathcal{D}_{\xi_2}\boldsymbol{X}_2\big) +
\big(\nwu{\mathcal{D}_{\xi_1}\boldsymbol{X}_1}\big) 
\big(\mathcal{D}_{\xi_2}\dot{\boldsymbol{X}}_2\big)  \\ - &
\big(\nwu{\mathcal{D}_{\xi_2}\dot{\boldsymbol{X}}_1}\big) 
\big(\mathcal{D}_{\xi_1}\boldsymbol{X}_2\big) -
\big(\nwu{\mathcal{D}_{\xi_2}\boldsymbol{X}_1}\big) 
\big(\mathcal{D}_{\xi_1}\dot{\boldsymbol{X}}_2\big)
\Big] ,
\end{aligned}
\end{equation}
where on the right hand side we have inserted the time derivative of the right hand side in (\ref{eq:J_appr}). 

Conversely, we always have the option to first specify the value of $\boldsymbol{\nabla}\cdot\dot{\boldsymbol{X}}$, for example by applying a set of physical space SBP operators as in (\ref{eq:DGCL}), and then use the relation (\ref{eq:GCL_disc}) constructively in order to define $\mathcal{J}$.

\subsection{Free-stream preservation}
\label{sec:curvi_FSP}
Again using the example of $2$D coordinates for brevity, SBP operators with respect to the physical coordinates in $\Omega$ can be defined as follows, only provided that $\mathcal{D}_{\xi_1}$ and $\mathcal{D}_{\xi_2}$ are SBP operators (\ref{eq:SBP_partial}) with respect to the reference domain quadrature $\hat{\mathcal{P}}$ ,
\begin{equation*}
\begin{aligned}
\mathcal{D}_{x_1} = & \frac{1}{2}\nwu{\mathcal{J}}^{-1}\Big[
\mathcal{D}_{\xi_1} \big(\nwu{\mathcal{D}_{\xi_2}\boldsymbol{X}_2}\big) +
\big(\nwu{\mathcal{D}_{\xi_2}\boldsymbol{X}_2}\big) \mathcal{D}_{\xi_1} -
\mathcal{D}_{\xi_2} \big(\nwu{\mathcal{D}_{\xi_1}\boldsymbol{X}_2}\big) -
\big(\nwu{\mathcal{D}_{\xi_1}\boldsymbol{X}_2}\big) \mathcal{D}_{\xi_2} 
\Big] \\
\mathcal{D}_{x_2}= & \frac{1}{2}\nwu{\mathcal{J}}^{-1}\Big[
\mathcal{D}_{\xi_2} \big(\nwu{\mathcal{D}_{\xi_1}\boldsymbol{X}_1}\big) +
\big(\nwu{\mathcal{D}_{\xi_1}\boldsymbol{X}_1}\big) \mathcal{D}_{\xi_2} -
\mathcal{D}_{\xi_1} \big(\nwu{\mathcal{D}_{\xi_2}\boldsymbol{X}_1}\big) -
\big(\nwu{\mathcal{D}_{\xi_2}\boldsymbol{X}_1}\big) \mathcal{D}_{\xi_1} 
\Big] ,
\end{aligned}
\end{equation*}
see also \cite{Lundquist18,Alund20_JCP} for more details. Now assume that $\mathcal{D}_{\xi_1}$ and $\mathcal{D}_{\xi_2}$ are consistent, i.e.
\begin{equation*}
\mathcal{D}_{\xi_1} \mathbb{1} = 0, \quad \mathcal{D}_{\xi_2} \mathbb{1} = 0,
\end{equation*}
and assume furthermore that $\mathcal{D}_{\xi_1}$ and $\mathcal{D}_{\xi_2}$ commute, i.e. $\mathcal{D}_{\xi_1}\mathcal{D}_{\xi_2}=\mathcal{D}_{\xi_2}\mathcal{D}_{\xi_1}$. Notably,  this commuting property is automatically satisfied if the operators have a tensor (Kronecker) product structure, i.e. $\mathcal{D}_{\xi_1}=D_{\xi_1}\otimes I_{\xi_2}$, $\mathcal{D}_{\xi_2}=I_{\xi_1}\otimes D_{\xi_2}$. From these two assumptions of consistency and commutativity, we can write
\begin{equation*}
\begin{aligned}
\mathcal{D}_{x_1} \mathbb{1} = & \frac{1}{2}\nwu{\mathcal{J}}^{-1}\Big[
\mathcal{D}_{\xi_1} \big(\mathcal{D}_{\xi_2}\boldsymbol{X}_2\big)  -
\mathcal{D}_{\xi_2} \big(\mathcal{D}_{\xi_1}\boldsymbol{X}_2\big) 
\Big] = 0 \\
\mathcal{D}_{x_2}\mathbb{1}= & \frac{1}{2}\nwu{\mathcal{J}}^{-1}\Big[ 
\mathcal{D}_{\xi_2} \big(\mathcal{D}_{\xi_1}\boldsymbol{X}_1\big) -
\mathcal{D}_{\xi_1} \big(\mathcal{D}_{\xi_2}\boldsymbol{X}_1\big) 
\Big] = 0,
\end{aligned}
\end{equation*}
i..e. the consistency assumption (\ref{eq:freestream_space}) used in Lemma \ref{thm:FSP_semidisc} is satisfied. However,  condition (\ref{eq:DGCL}) of the same lemma clearly leads to a different value of $\boldsymbol{\nabla}\cdot\dot{\boldsymbol{X}} $ than the one in (\ref{eq:J_FSP}).  Hence the approximation (\ref{eq:J_appr}) is not suitable if free-stream preservation is desirable (the same conclusion of course also applies if the exact values of $J$ are used).  Instead we must take $\boldsymbol{\nabla}\cdot\dot{\boldsymbol{X}} $ as given from (\ref{eq:DGCL}), and then make constructive use of the initial value problem (\ref{eq:GCL_disc}) in order to define $\mathcal{J}$.  According to proposition \ref{thm:FSP_semidisc}, this leads to semi-discrete free-stream preservation. Also recall that, for fully discrete free-stream preservation we need to solve the semi-coupled system (\ref{eq:system}) for both $\mathcal{J}$ and $\boldsymbol{U}$ simultaneously by applying the same time integration scheme to both equations.

\section*{Acknowledgments}
This work is based on research supported in part by the National Research Foundation of South Africa (Grant Numbers: 89916).  Tomas Lundquist was also partly funded through Vetenskapsr{\aa}det, Sweden grant agreement 2020-03642 VR and the Swedish e-Science Research Centre (SeRC). Jan Nordstr\"{o}m was supported by Vetenskapsr\r adet, Sweden grant 2018-05084 VR and SeRC.
The opinions, findings and conclusions or recommendations expressed is that of the authors alone, and the NRF accepts no liability whatsoever in this regard.

\bibliographystyle{amsplain}
\bibliography{References1}
\end{document}

%% file: ex_shared.tex
\usepackage{lipsum}
\usepackage{amsfonts}
\usepackage{graphicx}
\usepackage{epstopdf}
\usepackage{algorithmic}
\ifpdf
  \DeclareGraphicsExtensions{.eps,.pdf,.png,.jpg}
\else
  \DeclareGraphicsExtensions{.eps}
\fi

\usepackage{bbold}
\usepackage{mathtools}

\newcommand{\nwu}[2][1]{\mathrlap{\underline{\phantom{\ensuremath{\mathrm{#2}\mkern-#1mu}}}}#2\mkern#1mu}

\newsiamremark{remark}{Remark}
\newsiamremark{hypothesis}{Hypothesis}
\crefname{hypothesis}{Hypothesis}{Hypotheses}
\newsiamthm{claim}{Claim}

\headers{A method-of-lines framework for energy stable ALE methods}{T.Lundquist, A. Malan, and J.Nordstr\"{o}m}

\title{A method-of-lines framework for energy stable arbitrary Lagrangian-Eulerian methods
\thanks{Preprint submitted to SIAM Journal on Numerical Analysis.}
}

\author{Tomas Lundquist\thanks{Department of
  Mathematics, Computational Mathematics, Link\"{o}ping University, Sweden
  (\email{tomas.lundquist@liu.se}).}
\and Arnaud Malan\thanks{Industrial CFD Research Group, Department of Mechanical Engineering, University of Cape Town, Private Bag X3, South Africa
  (\email{arnaud.malan@uct.ac.za}).}
\and Jan Nordstr\"{o}m\thanks{Department of
Mathematics, Link\"{o}ping University,  SE-581 83 Link\"{o}ping, Sweden and Department of Mathematics and Applied Mathematics,  University of Johannesburg, P.O. Box 524, Auckland Park 2006, South Africa (\email{jan.nordstrom@liu.se}).}
}

\usepackage{amsopn}